\pdfoutput=1
\RequirePackage{ifpdf}
\ifpdf 
\documentclass[pdftex]{sigma}
\else
\documentclass{sigma}
\fi

\numberwithin{equation}{section}

\newtheorem{Theorem}{Theorem}[section]
\newtheorem*{Theorem*}{Theorem}
\newtheorem*{PartI}{Part~I}
\newtheorem*{PartII}{Part~II}
\newtheorem*{PartIII}{Part~III}
\newtheorem{Corollary}[Theorem]{Corollary}
\newtheorem{Lemma}[Theorem]{Lemma}
\newtheorem{Proposition}[Theorem]{Proposition}
\newtheorem*{Claim*}{Claim}

{\theoremstyle{definition}
\newtheorem{Definition}[Theorem]{Definition}

\newtheorem{Example}[Theorem]{Example}
\newtheorem{Remark}[Theorem]{Remark}
\newtheorem{Hypothesis}[Theorem]{Hypothesis}

\newtheorem*{Convention*}{Convention}
}

\usepackage{scrextend} 
\usepackage{hycolor}
\usepackage{xcolor}

\newcommand{\C}{{\mathbb{C}}}
\newcommand{\D}{{\mathbb{D}}}

\newcommand{\F}{{\mathbb{F}}}

\newcommand{\N}{{\mathbb{N}}}

\newcommand{\R}{{\mathbb{R}}}
\renewcommand{\SS}{{\mathbb{S}}}

\newcommand{\Z}{{\mathbb{Z}}}
\newcommand{\Aa}{{\mathcal{A}}} 

\newcommand{\Ee}{{\mathcal{E}}}
\newcommand{\Ff}{{\mathcal{F}}}

\newcommand{\Ll}{{\mathcal{L}}} 
\newcommand{\Mm}{{\mathcal{M}}} 

\newcommand{\Pp}{{\mathcal{P}}}

\newcommand{\Xx}{{\mathcal{X}}}

\newcommand{\id}{{\rm id}} 

\newcommand{\diag}{{\rm diag}} 
\newcommand{\grad}{\mathop{\mathrm{grad}}} 
\newcommand{\INDEX}{\mathop{\mathrm{index}}} 

\newcommand{\cgraph}[1]{\Gamma_{\kern-.5ex{}#1}} 
\renewcommand{\Re}{{\rm Re}} 

\newcommand{\Crit}{\operatorname{Crit}} 
\newcommand{\spec}{\mathrm{spec}\,} 
\newcommand{\norm}{{\rm norm}}

\newcommand{\eps}{{\varepsilon}}

\newcommand{\GL}{{\rm GL}}

\renewcommand{\tt}{{\mathfrak t}} 
\newcommand{\INNER}[2]{\left\langle #1, #2\right\rangle}

\def\Nablatop#1{\nabla^{#1}\kern-.5ex{}}

\def\NABLA#1{{\mathop{\nabla\kern-.5ex\lower1ex\hbox{$#1$}}}}
\def\Nabla#1{\nabla\kern-.5ex{}_{#1}}
\def\Tabla#1{\Tilde\nabla\kern-.5ex{}_{#1}}
\def\Babla#1{\widebar\nabla\kern-.5ex{}_{#1}}
\def\abs#1{\mathopen|#1\mathclose|}
\def\Abs#1{\left|#1\right|}
\def\norm#1{\mathopen\|#1\mathclose\|}
\def\Norm#1{\left\|#1\right\|}

\renewcommand{\Tilde}{\widetilde}

\newcommand{\p}{{\partial}}

\newcommand{\INTO}{\hookrightarrow} 

\newcommand{\1}{{{\mathchoice {\rm 1\mskip-4mu l} {\rm 1\mskip-4mu l}
{\rm 1\mskip-4.5mu l} {\rm 1\mskip-5mu l}}}}

\makeatletter
\renewcommand*\env@matrix[1][\arraystretch]{%
 \edef\arraystretch{#1}%
 \hskip -\arraycolsep
 \let\@ifnextchar\new@ifnextchar
 \array{*\c@MaxMatrixCols c}}
\makeatother

\makeatletter
\let\save@mathaccent\mathaccent
\newcommand*\if@single[3]{%
 \setbox0\hbox{${\mathaccent"0362{#1}}^H$}%
 \setbox2\hbox{${\mathaccent"0362{\kern0pt#1}}^H$}%
 \ifdim\ht0=\ht2 #3\else #2\fi
 }
\newcommand*\rel@kern[1]{\kern#1\dimexpr\macc@kerna}
\newcommand*\widebar[1]{\@ifnextchar^{{\wide@bar{#1}{0}}}{\wide@bar{#1}{1}}}
\newcommand*\wide@bar[2]{\if@single{#1}{\wide@bar@{#1}{#2}{1}}{\wide@bar@{#1}{#2}{2}}}
\newcommand*\wide@bar@[3]{%
 \begingroup
 \def\mathaccent##1##2{%
 \let\mathaccent\save@mathaccent
 \if#32 \let\macc@nucleus\first@char \fi
 \setbox\z@\hbox{$\macc@style{\macc@nucleus}_{}$}%
 \setbox\tw@\hbox{$\macc@style{\macc@nucleus}{}_{}$}%
 \dimen@\wd\tw@
 \advance\dimen@-\wd\z@
 \divide\dimen@ 3
 \@tempdima\wd\tw@
 \advance\@tempdima-\scriptspace
 \divide\@tempdima 10
 \advance\dimen@-\@tempdima
 \ifdim\dimen@>\z@ \dimen@0pt\fi
 \rel@kern{0.6}\kern-\dimen@
 \if#31
 \overline{\rel@kern{-0.6}\kern\dimen@\macc@nucleus\rel@kern{0.4}\kern\dimen@}%
 \advance\dimen@0.4\dimexpr\macc@kerna
 \let\final@kern#2%
 \ifdim\dimen@<\z@ \let\final@kern1\fi
 \if\final@kern1 \kern-\dimen@\fi
 \else
 \overline{\rel@kern{-0.6}\kern\dimen@#1}%
 \fi
 }%
 \macc@depth\@ne
 \let\math@bgroup\@empty \let\math@egroup\macc@set@skewchar
 \mathsurround\z@ \frozen@everymath{\mathgroup\macc@group\relax}%
 \macc@set@skewchar\relax
 \let\mathaccentV\macc@nested@a
 \if#31
 \macc@nested@a\relax111{#1}%
 \else
 \def\gobble@till@marker##1\endmarker{}%
 \futurelet\first@char\gobble@till@marker#1\endmarker
 \ifcat\noexpand\first@char A\else
 \def\first@char{}%
 \fi
 \macc@nested@a\relax111{\first@char}%
 \fi
 \endgroup
}
\makeatother

\def\XXint#1#2#3{{\setbox0=\hbox{$#1{#2#3}{\int}$}
 \vcenter{\hbox{$#2#3$}}\kern-.5\wd0}}

\long\def\symbolfootnote[#1]#2{\begingroup%
\def\thefootnote{\fnsymbol{footnote}}\footnote[#1]{#2}\endgroup}

\newcommand{\gf}[1]{#1} 

\begin{document}

\allowdisplaybreaks

\newcommand{\arXivNumber}{2510.20945}

\renewcommand{\thefootnote}{}

\renewcommand{\PaperNumber}{032}

\FirstPageHeading

\ShortArticleName{The Linearized Floer Equation in a Chart}

\ArticleName{The Linearized Floer Equation in a Chart\footnote{This paper is a~contribution to the Special Issue on Geometry and Dynamics in memory of Will Merry. The~full collection is available at \href{https://sigma-journal.com/Merry.html}{https://sigma-journal.com/Merry.html}}}

\Author{Urs FRAUENFELDER~$^{\rm a}$ and Joa WEBER~$^{\rm b}$}

\AuthorNameForHeading{U.~Frauenfelder and J.~Weber}

\Address{$^{\rm a)}$~Universit\"at Augsburg, Germany}
\EmailD{\mail{urs.frauenfelder@math.uni-augsburg.de}}

\Address{$^{\rm b)}$~Unicamp, Campinas, Brazil}
\EmailD{\mail{joa@math.uni-bielefeld.de}}

\ArticleDates{Received October 31, 2025, in final form March 11, 2026; Published online April 01, 2026}

\Abstract{In this article, we are considering the Hessian of the area functional in a non-Darboux chart. This does not seem to have been considered before and leads to an interesting new mathematical structure which we introduce in this article and refer to as \emph{almost extendable weak Hessian~field}. Our main result is a Fredholm theorem for Robbin--Salamon operators
associated to \emph{non-continuous} Hessians which we prove by taking advantage of this new structure.}


\Classification{57R58; 47A53; 46B70}

\renewcommand{\thefootnote}{\arabic{footnote}}
\setcounter{footnote}{0}

\section{Introduction}

In a Darboux chart, the Hessian of the area functional is the
constant operator $A_0={\rm i}\p_t$.
In this note, we look at the Hessian of the area functional in a
\emph{non-Darboux} chart.
To the best of our knowledge, this was never done before.
Here, in particular, an additional summand shows up
whose discontinuity challenges the previous methods
to obtain the Fredholm property
for the associated Robbin--Salamon operators.

We discover an interesting structure from the point of view of scale
geometry. We think that the structure is of crucial importance in
order to understand the structure of Floer theory in general.
Therefore, we give this structure a name and refer to it as
an \emph{almost extendable weak Hessian field}.
We then study the Robbin--Salamon operator
$\D=\p_s +A$
associated to a~connecting path in an almost extendable weak Hessian
field $A$.

Our main result is that this Robbin--Salamon operator is Fredholm.
The difficulty in proving this result is that
the Hessian is not necessarily continuous.
Therefore, the improvement, due to Rabier~\cite{Rabier:2004a}, of the classical
Robbin--Salamon theorem~\cite{robbin:1995a} is not necessarily
applicable to the situation at hand.
However, the new notion of almost extendability discovered in this article
allows us to decompose the Hessian in the sum of two operators
where one is still continuous while the other one is of lower order,
in symbols
$A=F+C$.
In contrast to the Hessian itself, the two summands are not necessarily
symmetric any more.
Luckily the theorem of Rabier can as well deal with non-symmetric
situations as long as continuity is guaranteed.
We check that the conditions of Rabier apply to the operator
associated to the first summand in our situation
$\F=\p_s +F$.
Hence $\F$ is a Fredholm operator.
Since the second summand $C$ is of lower order, we show that it
gives rise to a compact multiplication operator. This then proves our
main result, because the Fredholm property is invariant under compact
perturbation.

This article is part of our endeavor to understand the basic structure
behind Floer theory
in order to make Floer theory applicable to a broader class of
problems involving Hamiltonian delay equation
as explained in~\cite{Frauenfelder:2024c,Frauenfelder:2024d}.

\section{The area functional in a Darboux chart}

On $\C^n$, the standard symplectic form
$\omega=\sum_{j=1}^n{\rm d}x_j\wedge {\rm d}y_j$ has the primitive
\[
\lambda=\frac12 \sum_{j=1}^n( x_j {\rm d}y_j- y_j {\rm d}x_j).
\]
We are considering the area functional
\[
 \Aa_0\colon \ C^\infty\bigl(\SS^1,\C^n\bigr)\to\R
 ,\qquad
 u\mapsto \int_{\SS^1} u^*\lambda,
\]
where $\SS^1=\R/\Z$.
We write $u\in C^\infty\bigl(\SS^1,\C^n\bigr)$ as a Fourier series
\[
 u(t)=\sum_{k\in\Z} u_k e^{2\pi {\rm i} kt}
 ,\qquad
 \dot u(t)=\sum_{\ell\in\Z} 2\pi {\rm i} \ell u_\ell e^{2\pi {\rm i} \ell t} ,
\]
where $u_k\in \C^n$ for $k\in\Z$.
Expressing the area functional via the
Fourier decomposition,\footnote{Let $z=x+{\rm i}y$ and $\xi=\hat x+{\rm i}\hat y$ where $x,y,\hat x,\hat y\in\R^n$.
 Hence $z=(x_1+{\rm i}y_1, \dots, x_n+{\rm i}y_n)$ and
 $\xi=(\hat x_1+{\rm i}\hat y_1,\allowbreak\dots, \hat x_n+{\rm i}\hat y_n)$ and we calculate
 \smash{$
 \Re\INNER{z}{{\rm i}\xi}_{\C^n}
 =\Re\sum_{k=1}^n(x_k+{\rm i}y_k)\overline{(-\hat y_k+{\rm i}\hat x_k)}
 =-\lambda_z\xi .
 $}
 }
\begin{equation*}
\begin{split}
 \Aa_0(u)
 &=\int_0^1 \lambda_u\dot u\, {\rm d}t
\\
 &=-\int_0^1\frac12 \Re \INNER{u}{{\rm i} \dot u}_{\C^n} {\rm d}t
\\
 &=-\int_0^1\frac12 \Re \sum_{k,\ell\in\Z} 2\pi\ell
 \big\langle u_k e^{2\pi {\rm i} kt},{\rm i}^2 u_\ell e^{2\pi {\rm i} \ell t}\big\rangle_{\C^n} {\rm d}t
\\
 &=\Re \sum_{k,\ell\in\Z} \pi\ell
 \INNER{u_k}{u_\ell}_{\C^n}
 \underbrace{\int_0^1 e^{2\pi {\rm i} (k-\ell) t} \, {\rm d}t}_{=\delta_{k\ell}}
\\
 &=\sum_{k\in\Z} \pi k \abs{u_k}^2_{C^n} .
\end{split}
\end{equation*}
In particular, we see that the area functional in a Darboux chart is just
a quadratic function in the Fourier coefficients.
Therefore, the Hessian is constant independent of $u$
with eigenvalues~$2\pi k$ where $k\in\Z$.
Abbreviating \smash{$H_k=W^{k,2}\bigl(\SS^1,\C^n\bigr)$},
we see that the Hessian of $\Aa_0$ at every point
$u\in H_1$ is a Fredholm operator of index zero
from $H_k$ to $H_{k-1}$ for every $k\in\N$.

\section{Euclidean inner product and symplectic forms}

\subsection{Associated anti-symmetric matrix}

\begin{Definition}
Let $\INNER{\cdot}{\cdot}$ be the Euclidean inner product on $\R^{2n}$.
Let $\mathfrak{U}\subset\R^{2n}$ be an open subset carrying an exact symplectic
form $\omega={\rm d}\lambda$. Then the identity
\begin{equation}\label{eq:B}
 \INNER{\cdot}{\cdot}
 =\omega_x(\cdot,B_x\cdot)
\end{equation}
for $x\in \mathfrak{U}$ determines a map $\mathfrak{U}\to\Ll\bigl(\R^{2n}\bigr)$, $x\mapsto B_x$. Each
linear map $B_x=B_x^\omega$ is unique by non-degeneracy of~$\omega_x$.
\end{Definition}

\begin{Lemma}\label{le:hghjghj8656h}
At any point $x\in \mathfrak{U}$ the linear map
$B_x\colon\R^{2n}\to\R^{2n}$ is invertible and satisfies
\begin{itemize}\itemsep=0ex
\item[$(i)$]
 $\omega(\xi,B_x\eta)=-\omega(B_x\xi,\eta)$
$(\omega$-anti-symmetry of $B_x)$,
\item[$(ii)$]
 $\INNER{\xi}{B_x\eta}=-\INNER{B_x\xi}{\eta}$
 $($anti-symmetry of $B_x)$,
\item[$(iii)$]
 $\INNER{\xi}{-B_x^2\eta}=\INNER{-B_x^2\xi}{\eta}$
 $($symmetry of $-B_x^2)$,
\item[$(iv)$]
 $\INNER{-B_x^2\xi}{\xi}=\abs{B\xi}^2>0$, $\xi\not= 0$
 $($positive definiteness of $-B_x^2)$,
\item[$(v)$]
 $\sqrt{-B_x^2}$ and its inverse are both symmetric positive
 definite
\end{itemize}
for all $\xi,\eta\in\R^{2n}$.
\end{Lemma}

As we shall see in Remark~\ref{rem:B+} further below,
the invertible linear maps $B_x\colon\R^{2n}\to\R^{2n}$ are
orientation preserving, in symbols $B_x\in\GL^+(2n,\R)$.

\begin{proof}
Suppose $B_x\eta=0$, then it follows from~(\ref{eq:B}) and
non-degeneracy of $\INNER{\cdot}{\cdot}$ that $\eta=0$.
Hence $B_x\colon\R^{2n}\to\R^{2n}$ is injective and therefore surjective.

(i)~By~(\ref{eq:B}) and symmetry of $\INNER{\cdot}{\cdot}$
and anti-symmetry of $\omega$, we obtain
\begin{equation*}
\begin{split}
 \omega_x(\xi,B_x\eta)
 =\INNER{\xi}{\eta}
 =\INNER{\eta}{\xi}
 =\omega_x(\eta,B_x\xi)
 =-\omega_x(B_x\xi,\eta) .
\end{split}
\end{equation*}

(ii)~By~(\ref{eq:B}) and anti-symmetry of $\omega$ and symmetry of
$\INNER{\cdot}{\cdot}$, we obtain
\begin{equation*}
\begin{split}
 \INNER{\xi}{B_x\eta}
 &=\omega_x\bigl(\xi,B_x^2\eta\bigr)
 \stackrel{\text{(i)}}{=}\omega_x\bigl((-B_x)^2\xi,\eta\bigr)
 =-\omega_x\bigl(\eta,B_x^2\xi\bigr)
 =-\INNER{\eta}{B_x\xi}
 =\INNER{-B_x\xi}{\eta} .
\end{split}
\end{equation*}

(iii)~holds by applying twice~(ii).

(iv)~By~(ii) and injectivity, we get
\[
 \INNER{-B_x^2\xi}{\xi}
 =\INNER{B_x\xi}{B_x\xi}
 =\abs{B_x\xi}^2>0.
\]

(v)~By Heron's construction of the square root of a symmetric positive
definite matrix $Q$, Lemma~\ref{le:Heron}, the square root is
symmetric positive definite and depends smoothly on $Q$.
Taking the inverse preserves symmetry and positive definiteness.
\end{proof}

For the previous and the following lemma,
see also~\cite[Proposition~2.5.6]{McDuff:2017b}.

\subsection{Compatible almost complex structure}

\begin{Lemma}\label{le:J_B}
The linear map \smash{$J_B:=\sqrt{-B^2}^{\,-1} B\colon\R^{2n}\to\R^{2n}$}
is an almost complex structure on~$\mathfrak{U}$ compatible with $\omega$, in symbols
\begin{itemize}\itemsep=0ex
\item[$(a)$]
 $J_BJ_B=-\1$
 $($almost complex structure$)$,
\item[$(b)$]
 $g_{J_B}(\xi,\eta):=\omega(\xi,J_B \eta)=g_{J_B}(\eta,\xi)$ $(g_{J_B}$ symmetric$)$,
\item[$(c)$]
 $g_{J_B}(\xi,\xi)=\big\langle{\xi},{\sqrt{-B^2}^{\,-1} \xi}\big\rangle>0$, $\xi\not=0$
 $(g_{J_B}$ positive definite$)$
\end{itemize}
for all vector fields $\xi$, $\eta$ along $\mathfrak{U}\subset\R^{2n}$.
\end{Lemma}

\begin{proof}
(a) Since $B$ commutes with $-B^2$, then by Corollary~\ref{cor:Heron}
it also commutes with $\sqrt{-B^2}$, hence with the inverse of \smash{$\sqrt{-B^2}$}.
This is used in equality 2 below
\begin{align*}
 J_BJ_B
 &=\sqrt{-B^2}^{\,-1} B \sqrt{-B^2}^{\,-1} B
 \stackrel{2}{=}\sqrt{-B^2}^{\,-1}\sqrt{-B^2}^{\,-1} BB\\
 &=\bigl(\sqrt{-B^2}\sqrt{-B^2}\bigr)^{-1} BB
 =(-BB)^{-1} BB=-\1 .
\end{align*}

(b) Equality 3 below holds analogously to equality 2 in (a),
equality 4 uses~(\ref{eq:B}), equality 5 is by
Lemma~\ref{le:hghjghj8656h} part~(v),
equality 6 is by symmetry of the Euclidean inner product:{\samepage
\begin{align*}
\begin{split}
 \underline{g_{J_B}(\xi,\eta)}
 :={}&\omega(\xi,J_B\eta)
 =\omega\big(\xi, \sqrt{-B^2}^{\,-1} B\eta\big)
 \stackrel{3}{=}\omega\big(\xi,B\sqrt{-B^2}^{\,-1}\eta\big)
 \stackrel{4}{=}\underline{\big\langle{\xi},{\sqrt{-B^2}^{\,-1}\eta}\big\rangle}\\
 \stackrel{5}{=}{}&\big\langle{\sqrt{-B^2}^{\,-1}\xi},{\eta}\big\rangle
 =\big\langle{\eta},{\sqrt{-B^2}^{\,-1}\xi}\big\rangle
 =g_{J_B}(\eta,\xi) .
 \end{split}
\end{align*}
The last step uses equality of the two underlined terms,
just commute $\xi$ and $\eta$.}

(c) Step 2 below holds analogously to equality 2 in (a),
Step~3 by~(\ref{eq:B}), namely
\begin{equation*}
 g_{J_B}(\xi,\xi)
 :=\omega\big(\xi, \sqrt{-B^2}^{\,-1} B\xi\big)
 \stackrel{2}{=}\omega\big(\xi, B\sqrt{-B^2}^{\,-1} \xi\big)
 \stackrel{3}{=}\big\langle{\xi},{\sqrt{-B^2}^{\,-1}\xi}\big\rangle
 >0
\end{equation*}
pointwise at each $x\in \mathfrak{U}$ whenever $\xi(x)\not=0$.
This proves Lemma~\ref{le:J_B}.
\end{proof}

\begin{Lemma}\label{le:J-det=1}
Any complex structure\footnote{A linear map $J\colon V\to V$ such that $J^2=-\1$.}
$J$ on a real vector space $V$ of finite dimension $2n$ has determinant $1$.
\end{Lemma}

\begin{proof}
For a complex basis on $V$
the matrix representation of $J$ is the complex
$n\times n$ matrix $J_\C={\rm i}\1_n$.
Its determinant is $\det (J_\C)={\rm i}^n$.
In general, writing a complex $n\times n$ matrix $Z=X+{\rm i}Y$
as sum of real and imaginary parts, define the corresponding real
$2n\times 2n$ matrix as follows and observe the identity
\[
 Z_\R:=\begin{pmatrix}X&-Y\\Y&X\end{pmatrix}
 =\begin{pmatrix}\1&\hphantom{-}\1\\ {\rm i}\1&-{\rm i}\1\end{pmatrix}^{-1}
 \begin{pmatrix}X-{\rm i}Y&0\\0&X+{\rm i}Y\end{pmatrix}
 \begin{pmatrix}\1&\hphantom{-}\1\\ {\rm i}\1&-{\rm i}\1\end{pmatrix}.
\]
Take the determinant to obtain the formula
\[
 \det Z_\R
 =\det \begin{pmatrix}\bar Z&0\\0&Z\end{pmatrix}
 =\det\bar Z\cdot \det Z
 =\overline{\det Z}\cdot \det Z
 =\abs{\det Z}^2 .
\]
Hence $\det J=\abs{\det J_\C}^2=\abs{{\rm i}^n}^2=1$.
\end{proof}

\begin{Remark}\label{rem:B+}
The linear map $B_x$ determined by~(\ref{eq:B})
has positive determinant
\begin{equation}\label{eq:B+}
 \det B
 =\det\bigl(\sqrt{-B^2} J_B\bigr)
 =\det\sqrt{-B^2} \cdot\det J_B
 =\det\sqrt{-B^2}
 > 0 .
\end{equation}
Here identity one is by definition of $J_B$ and identity three by
Lemma~\ref{le:J-det=1}.
Strictly positive holds true by
Lemma~\ref{le:hghjghj8656h}\,(v).
\end{Remark}

\section{Action functional}
\label{sec:unp-action}

We denote by $\SS^1$ the circle $\R/\Z$.
Consider the Hilbert space triple
\begin{equation}\label{eq:triple}
 H_0:=L^2\bigl(\SS^1,\R^{2n}\bigr),\qquad
 H_1:=W^{1,2}\bigl(\SS^1,\R^{2n}\bigr),\qquad
 H_2:=W^{2,2}\bigl(\SS^1,\R^{2n}\bigr).
\end{equation}
Given an open subset $\mathfrak{U}\subset\R^{2n}$,
let $\omega={\rm d}\lambda$ be an exact symplectic form on $\mathfrak{U}$.
Define open subsets of $H_1$ and $H_2$ by
\begin{equation}\label{eq:U_2}
 U_\ell:=\bigl\{u\in H_\ell\mid u(t)\in\mathfrak{U}\; \forall t\in\SS^1\bigr\}
 \subset C^0\bigl(\SS^1,\mathfrak{U}\bigr),\qquad \ell=1,2.
\end{equation}
\begin{Convention*}
We write
\smash{$\dot u(t):=\tfrac{{\rm d}}{{\rm d}t} u(t)$}.
All sums in Section~\ref{sec:unp-action} run from $1$ to $2n$.
\end{Convention*}

\begin{Definition}
The \textit{symplectic action} (functional) is defined by
\begin{equation}\label{le:action_0}
 \Aa_0\colon\ H_1\supset U_1\to\R
 ,\qquad
 u\mapsto \int_{\SS^1} u^*\lambda .
\end{equation}
\end{Definition}

\subsection{First derivative and gradient}

\begin{Lemma}[gradient]\label{le:gradient}
The derivative at $u\in U_1$ in direction $\xi\in H_1$ is
\[
 {\rm d}\Aa_0|_u\xi
 =\int_0^1\omega_{u_t}(\xi_t,\dot u_t)\, {\rm d}t
 =\int_0^1\big\langle{\xi_t},{B^{-1}_{u_t}\dot u_t}\big\rangle {\rm d}t
 =:\big\langle{\xi},{B^{-1}_{u}\dot u}\big\rangle_{H_0},
\]
where $B_{u_t}\in\Ll\bigl(\R^{2n}\bigr)$ is determined by~\eqref{eq:B}.
So the \textit{$\bigl(L^2\bigr)$-gradient} is given by
\[
 (\grad \Aa_0|_u )(t)=B_{u(t)}^{-1}\dot u(t)
\]
pointwise for every $t\in\SS^1$.
\end{Lemma}

\begin{proof}[Proof (global version)]
Given $u\in U_1$ and $\xi\in H_1$,
let $u_r$ be a smooth family with $u_0=u$ and
\smash{$\tfrac{{\rm d}}{{\rm d}r}\big|_{r=0}u_r=\xi$}. We compute
\begin{align*}
 {\rm d}\Aa_0|_u\xi
 &=\frac{{\rm d}}{{\rm d}r}\Big|_{r=0}\Aa_0(u_r)
 =\int_{\SS^1}\frac{{\rm d}}{{\rm d}r}\Big|_{r=0} u_r^*\lambda
 \stackrel{3}{=}\int_{\SS^1} u^* L_\xi\lambda\\
 &=\int_{\SS^1} u^* ({\rm d}{i}_\xi+i_\xi {\rm d})\lambda
 \stackrel{5}{=}\int_{\SS^1} {\rm d}u^*i_\xi\lambda +\int_{\SS^1} u^*i_\xi \omega
 \stackrel{6}{=} 0+\int_0^1\omega_{u_t}(\xi_t,\dot u_t)\, {\rm d}t .
\end{align*}
Here Step~3 is by definition of the Lie derivative
and Step~4 is Cartan's formula.
Step~5 uses that the exterior derivative ${\rm d}$ commutes with pull-back.
Step~6 is by Stokes' theorem and the fact that the integral is over the
empty set $\p\SS^1=\varnothing$, so the integral is $0$.
This proves Lemma~\ref{le:gradient}.
\end{proof}

\begin{proof}[Proof (local version)]
Let $u\colon\SS^1\to \mathfrak{U}$ be in $U_1$ and
$\xi\colon\SS^1\to\R^{2n}$ in $H_1$.
Take the derivative of $t\mapsto (\lambda|_u\xi)(t)$ to get
\begin{equation*}
\begin{split}
 \frac{{\rm d}}{{\rm d}t} ((u^*\lambda) \xi )(t)
 &=\frac{{\rm d}}{{\rm d}t}\sum_i\lambda_i(u(t))\xi_i(t) \\
 &=\sum_{i,j} \p_j\lambda_i(u(t))\dot u_j(t) \xi_i(t)\, {\rm d}t
 +\sum_i \lambda_i(u(t)) \dot\xi_i(t)\, {\rm d}t .
\\
 &\stackrel{3}{=}\sum_{i,j} \p_i\lambda_j(u(t))\dot u_i(t) \xi_j(t)\, {\rm d}t
 +\underline{\sum_i \lambda_i(u(t)) \dot\xi_i(t)\, {\rm d}t},
\end{split}
\end{equation*}
where in Step 3 we renamed the summation indices $(i,j)$ by $(j,i)$.
Let $u_r$ be a variation associated to $u$ and $\xi$.
Substitute the underlined term in what follows
\begin{align}
 \frac{{\rm d}}{{\rm d}r}\Big|_{r=0}
 (u_r^*\lambda )(t)
 &\stackrel{1}{=}\sum_{i,j} \p_j\lambda_i(u(t)) \xi_j(t) \dot u_i(t)\, {\rm d}t
 +\underline{\sum_i \lambda_i(u(t)) \dot\xi_i(t)\, {\rm d}t}\nonumber\\
 &=\sum_{i,j}
 \Bigl(
 \underbrace{\p_j\lambda_i(u(t)) \xi_j(t) \dot u_i(t)
 -\p_i\lambda_j(u(t))\dot u_i(t) \xi_j(t)}
 _{=\omega_{ji}(u(t)) \xi_j(t)\dot u_i(t)}
 \Bigr) {\rm d}t
 +\tfrac{{\rm d}}{{\rm d}t} ((u^*\lambda) \xi )(t) .\label{eq:hjvghv46yg}
\end{align}
We integrate and use $1$-periodicity in time $t$ to obtain the formula
\begin{equation*}
\begin{split}
 {\rm d}\Aa_0|_u\xi
 =\int_{\SS^1}\frac{{\rm d}}{{\rm d}r}\Big|_{r=0} u_r^*\lambda
 &=\int_0^1 \sum_{i,j} \omega_{ji}(u_t) \xi_j(t)\dot u_i(t) \, {\rm d}t
 =\int_0^1\omega_{u_t}(\xi_t,\dot u_t)\, {\rm d}t .
\end{split}
\end{equation*}
This proves Lemma~\ref{le:gradient}.
\end{proof}

\subsection{Second derivative}

Consider an open subset $\mathfrak{U}$
of $\R^{2n}$ with coordinates $x=(x_1,\dots,x_{2n})$.
A primitive $\lambda$ of $\omega$
is of the form \smash{$\lambda=\sum_{i=1}^{2n} \lambda_i(x) {\rm d}x_i$}.
We denote the second derivatives of the coefficient functions
$x\mapsto \lambda_i(x)$ and their differences by
\[
 \Lambda_{kji}(x):=\p_k\p_j\lambda_i(x)
 ,\qquad
 L_{kji}(x):=\Lambda_{kji}(x)-\Lambda_{ikj}(x)
 ,\qquad
 i,j,k=1,\dots,2n,
\]
at any point $x\in\mathfrak{U}\subset\R^{2n}$.

To compute the Hessian of the function
$\Aa_0\colon H_1\supset U_1\to\R$
we introduce the following covariant tensor field of type $(0,3)$
in the notation of~\cite[Section~2]{oneill:1983a}.
The set of vector fields, respectively, functions, along
$\mathfrak{U}\subset\R^{2n}$ are denoted by
\[
 \Xx(\mathfrak{U}):=C^\infty\bigl(\mathfrak{U},\R^{2n}\bigr)
 ,\qquad
 \Ff(\mathfrak{U}):=C^\infty(\mathfrak{U},\R) .
\]

\begin{Definition}
A tensor field of type $(0,3)$ is defined by the sum
\begin{equation*}
\begin{split}
 L\colon\ \Xx(\mathfrak{U}) \times \Xx(\mathfrak{U})\times
 \Xx(\mathfrak{U})
 &\to \Ff(\mathfrak{U}),
 \qquad
 (\eta,\xi,\zeta)
 \mapsto
 \sum_{i,j,k} L_{kji} \eta_k\xi_j \zeta_i .
\end{split}
\end{equation*}
\end{Definition}

\begin{Remark}[Darboux charts]
Locally, in Darboux coordinates $(q_1,\dots,q_n,p_1,\dots, p_n)$,
the symplectic form $\omega$ has as a primitive the Liouville form
$\lambda=\sum_{i=1}^n p_i {\rm d}q_i$.
Observe that all coefficients, namely the $p_i$, are linear functions
in the coordinates. Therefore, the tensor field vanishes $L\equiv 0$.
\end{Remark}

We discuss some symmetries of $L$.
By the theorem of Schwarz,
$
 \Lambda_{kji}=\Lambda_{jki}
$
and therefore
\begin{equation}\label{eq:Schwarz}
 L_{kji}-L_{jki}-L_{ijk}=0 .
\end{equation}
Furthermore, by definition of the $L_{kji}$, the sum of cyclic
permutations vanishes
\[
 L_{kji}+L_{ikj}+L_{jik}=0 .
\]
Consequently, the sum of cyclic permutations in $L$ vanishes as well
\begin{equation*}
 L(\eta,\xi,\zeta)+L(\zeta,\eta,\xi)
 +L(\xi,\zeta,\eta)
 =0 .
\end{equation*}
By the consequence~(\ref{eq:Schwarz}) of Schwarz,
the following sum vanishes as well:
\begin{equation*}
 L(\eta,\xi,\zeta)-L(\xi,\eta,\zeta)
 -L(\zeta,\xi,\eta)
 =0 .
\end{equation*}

\begin{Lemma}\label{le:2nd-deriv}
The second derivative of $\Aa_0\colon U_1\to\R$
at $u$ evaluated on two vector fields
$\xi,\eta\in H_1$ is the symmetric bilinear form given by
\begin{equation*}
\begin{split}
 {\rm d}^2\Aa_0|_u\left(\xi,\eta\right)
 = 
 \int_0^1 \omega_{u}(\xi,\dot\eta)\, {\rm d}t
 +\int_0^1 L_u(\eta,\xi,\dot u)\, {\rm d}t .
\end{split}
\end{equation*}
\end{Lemma}

\begin{proof}
Given $u\in U_1$ and $\xi,\eta\in H_1$,
let $u_{r,\rho}$ be an associated two-parameter variation, that is,
$u_{0,0}=u$ and
\[
 \frac{{\rm d}}{{\rm d}r} u_{r,\rho}(t)\Big|_{(r,\rho)=(0,0)}=\xi(t)
 ,\qquad
 \frac{{\rm d}}{{\rm d}\rho} u_{r,\rho}(t)\Big|_{(r,\rho)=(0,0)}=\eta(t) .
\]
We use equality 1 in~(\ref{eq:hjvghv46yg}) in Step~2 to compute
{\samepage\begin{align}
 \frac{{\rm d}}{{\rm d}\rho}\frac{{\rm d}}{{\rm d}r}\Big|_{(r,\rho)=(0,0)}
 \bigl(u_{r,\rho}^*\lambda \bigr)
 &= \frac{{\rm d}}{{\rm d}\rho}\Big|_{\rho=0}\biggl(
 \sum_{i,j} \p_j\lambda_i\bigl(u_{(0,\rho)}\bigr) \xi_j \bigl(\dot u_{(0,\rho)}\bigr)_i\, {\rm d}t
 +\sum_i \lambda_i\bigl(u_{(0,\rho)}\bigr) \dot\xi_i\, {\rm d}t\biggr) \nonumber\\
 &=\sum_{i,j,k}\p_k\p_j\lambda_i(u) \eta_k \xi_j \dot u_i\, {\rm d}t
 +\sum_{i,j}\p_j\lambda_i(u)
 \bigl(\xi_j\dot\eta_i+\eta_j\dot\xi_i\bigr) {\rm d}t\label{eq:jhghughu77546}
\end{align}
evaluated pointwise at $t$.}
In the following calculation,
we use the definition~(\ref{le:action_0}) of the action functional
in Step~2, and~(\ref{eq:jhghughu77546}) in Step~3, to obtain
\begin{equation*}
\begin{split}
 {\rm d}^2\Aa_0(u)\bigl(\xi,\eta\bigr)
 ={}&\frac{{\rm d}}{{\rm d}\rho}\frac{{\rm d}}{{\rm d}r}\Big|_{(r,\rho)=(0,0)}\Aa_0(u_{r,\rho})\\
 \stackrel{2}{=}{}&
 \int_{\SS^1} \frac{{\rm d}}{{\rm d}\rho}\frac{{\rm d}}{{\rm d}r}\Big|_{(r,\rho)=(0,0)}
 {u_{r,\rho}}^*\lambda \\
 \stackrel{3}{=}{}&
 \int_0^1 \biggl(\sum_{i,j,k}\p_k\p_j\lambda_i(u) \eta_k \xi_j \dot u_i\, {\rm d}t
 +\sum_{i,j}\p_j\lambda_i(u)
 \bigl(\xi_j\dot\eta_i+\eta_j\dot\xi_i\bigr)\biggr) {\rm d}t
\\
 \stackrel{4}{=}{}&
 \int_0^1 \biggl(\sum_{i,j,k}
 \underbrace{\p_k\p_j\lambda_i(u)}_{=:\Lambda_{kji}(u)}
 \eta_k \xi_j \dot u_i\, {\rm d}t
 +\sum_{i,j}
 \underbrace{\p_j\lambda_i(u)
 \bigl(\xi_j\dot\eta_i-\dot\eta_j\xi_i\bigr)}
 _{=\omega_{ji}(u) \xi_j\dot \eta_i}
 \biggr) {\rm d}t
 \\
 &{}
 {-}\,\int_0^1 \sum_{i,j,k}\p_k\p_j\lambda_i(u) \dot u_k \eta_j \xi_i\, {\rm d}t
 \\
 ={}& 
 \int_0^1 \omega_{u}(\xi,\dot\eta)\, {\rm d}t
 +\int_0^1 \sum_{i,j,k}
 (\underbrace{\Lambda_{kji}(u)-\Lambda_{ikj}(u)}_{=:L_{kji}(u)})
 \eta_k \xi_j \dot u_i
 \, {\rm d}t,
\end{split}
\end{equation*}
where Step~4 is by integration by parts
(boundary terms vanish by periodicity).
Symmetry holds by the theorem of Schwarz.
This proves Lemma~\ref{le:2nd-deriv}.
\end{proof}

\subsection{Para-Darboux Hessian}

Using the metric isomorphism, we turn the $(0,3)$ tensor $L$
into the $(1,2)$ tensor
\begin{equation}\label{eq:L-(1,2)}
\begin{split}
 \bar L\colon\ \Xx(\mathfrak{U}) \times \Xx(\mathfrak{U})
 &\to \Xx(\mathfrak{U})
 ,\qquad
 (\eta,\zeta)\mapsto \bar L(\eta,\zeta)
\end{split}
\end{equation}
determined by the identity
\[
 \INNER{\bar L(\eta,\zeta)}{\xi}=L(\eta,\xi,\zeta) .
\]
Now we can write the Hessian of $\Aa_0$ at $u\in U_1$ and for
vector fields $\xi,\eta\in H_1$ in terms of the $L^2$ inner product.
Namely, by Lemma~\ref{le:2nd-deriv} we get
\begin{equation}\label{eq:Hessian-A_0}
 {\rm d}^2\Aa_0|_u (\xi,\eta )
 =\INNER{\xi}{B^{-1}_u\dot\eta +\bar L_u(\eta,\dot u)}_{H_0}
 =\INNER{\xi}{A^u_0\eta}_{H_0},
\end{equation}
where $H_0=L^2\bigl(\SS^1,\R^{2n}\bigr)$ and the loop
$t\mapsto B_{u(t)}\in\GL^+\bigl(2n,\R^{2n}\bigr)$ is
determined pointwise at $t$ by~(\ref{eq:B}) and the last identity
determines the linear operator~$A^u_0$.

\begin{Definition}
The \textit{para-Darboux Hessian operator} of the action $\Aa_0$
at $u\in U_1$ is, due to~(\ref{eq:Hessian-A_0}), the bounded
linear map
\begin{equation}\label{eq:Hess-op-0}
\begin{split}
 A_0^u\colon\ H_1\to H_0 ,\qquad
 \eta\mapsto B_u^{-1}\dot \eta+\bar L_u(\eta,\dot u) .
\end{split}
\end{equation}
For any $u\in U_2$, see~(\ref{eq:U_2}), this is a bounded linear map
$A_0^u|_{H_2}\colon H_2\to H_1$.
\end{Definition}

\begin{Lemma}[$H_0$-symmetry]\label{le:A^u_0-symm}
$\forall u\in U_1\colon
 \INNER{\xi}{A^u_0\eta}_{H_0}=\INNER{A^u_0\xi}{\eta}_{H_0}$
$\forall \xi,\eta\in H_1$.
\end{Lemma}

\begin{proof}
By Lemma~\ref{le:2nd-deriv}, the Hessian is symmetric, now use the
identity~(\ref{eq:Hessian-A_0}).
\end{proof}

\begin{Remark}[Morse--Bott]\label{rem:Morse--Bott}
The critical points of $\Aa_0$ are precisely the constant loops, in symbols
\[
 \Crit\Aa_0=\mathfrak{U}.
\]
Let $u$ be a constant loop in $\mathfrak{U}$.
Then the second term of the para-Darboux Hessian operator vanishes.
Hence the para-Darboux Hessian operator is just given by
\[
A_u^0\eta=B_u^{-1}\dot\eta \qquad \text{and} \qquad
 \ker A_0^u=\{\eta\in H_1\mid \dot \eta=0\}
 \simeq \R^{2n}
 =T_u\mathfrak{U}
.\]
\end{Remark}

\section{Perturbed action functional}

As discussed in Remark~\ref{rem:Morse--Bott},
the functional $\Aa_0$ is not Morse, but only Morse--Bott.
In order to get a Morse functional, we look at perturbations
of the area functional $\Aa_0$.
To this end, we introduce time-dependent Hamiltonians
which are $1$-periodic in time.
The Hilbert spaces $(H_0,H_1,H_2)$ are defined by~(\ref{eq:triple})
and open subsets $U_1\subset H_1$ and $U_2\subset H_2$ by~(\ref{eq:U_2}).

\begin{Definition}[perturbed action]
For $C^3$ functions $h\colon\SS^1\times\mathfrak{U}\to\R$,
notation $h_t(x):=h(t,x)$,
the \textit{perturbed action functional} is defined by
\[
 \Aa_h\colon\ U_1\to\R
 ,\qquad
 u\mapsto\int_{\SS^1}u^*\lambda-\int_0^1 h_t(u(t))\, {\rm d}t .
\]
\end{Definition}

To define the \textit{Hamiltonian vector field} $X_{h_t}$, we choose the
convention that
\[
 {\rm d}h_t(\cdot)=\omega(\cdot,X_{h_t})
\]
whenever $t\in\SS^1$.
The derivative of $\Aa_h$ at $u\in U_1$ in direction $\xi\in H_1$ is
given by
\[
 {\rm d}\Aa_h|_u\xi
 =\int_0^1 (\omega_{u_t}(\xi_t,\dot u_t)
 -{\rm d}h_t|_{u_t}\xi_t )\, {\rm d}t
 =\int_0^1\omega_{u_t}(\xi_t,\dot u_t-X_{h_t}(u_t))\, {\rm d}t .
\]
In particular, critical points of $\Aa_h$ are
solutions of the ODE
\[
 \dot u(t)=X_{h_t}(u(t))
\]
for $t\in\SS^1$, i.e., $1$-periodic orbits of the Hamiltonian vector
field of $h$.

\subsection{Perturbed para-Darboux Hessian}

Given $t\in\SS^1$, the \textit{Hessian operator of the function}
$h_t\colon\mathfrak{U}\to\R$ at $x\in\mathfrak{U}$ is the linear map
whose matrix with respect to the canonical basis is given by
\begin{equation}\label{eq:matrix-a_t}
 a_t|_x=(\p_i\p_j h_t|_x)_{i,j=1}^{2n}\colon\ \R^{2n}\to\R^{2n} .
\end{equation}
It is determined by the identity
\[
 {\rm d}^2h_t|_x(v,w)
 =\INNER{v}{a_t|_x w}
 ,\qquad
 {\rm d}^2h_t|_x(v,w)
 :=\frac{{\rm d}}{{\rm d}\tau}\Big|_{0}\frac{{\rm d}}{{\rm d}\eps}\Big|_{0}
 h_t(x+\eps v+\tau w),
\]
for all $v,w\in\R^{2n}$. By Schwarz's theorem, ${\rm d}^2h_t|_x$ is
symmetric, so $a_t|_x={a_t|_x}^T$.

\begin{Definition}
The \textit{Hessian operator of the Hamiltonian perturbation} $h$
at $u\in U_1$ is the bounded linear map $a^u\colon H_0\to H_0$
defined pointwise at $t$ by
\begin{equation}\label{eq:a_u}
 (a^u\eta )(t)=a_t|_{u_t}\eta_t .
\end{equation}
Still for $u\in U_1$,
this map is also bounded as a map $H_1\to H_1$.
Since the matrix is symmetric, it holds that
$ \INNER{\xi}{a^u_0\eta}_{H_0}=\INNER{a^u_0\xi}{\eta}_{H_0}$
for all vector fields $\xi,\eta\in H_0$.
\end{Definition}

The second derivative of the perturbed action at $u\in U_1$ is, as
a consequence of Lemma~\ref{le:2nd-deriv},
for $\xi,\eta\in H_1$ given by
\begin{equation*}
\begin{split}
 {\rm d}^2\Aa_h(u) (\xi,\eta )
 =\int_0^1\bigl( \omega_{u_t}(\xi_t,\dot\eta_t)
 + L_{u_t}(\eta_t,\xi_t,\dot u_t)
 - {\rm d}^2h|_{u_t}(\eta_t,\xi_t)
 \bigr) {\rm d}t .
\end{split}
\end{equation*}

\begin{Definition}
The \textit{perturbed para-Darboux Hessian operator}
at $u\in U_1$ is the bounded linear map given by the difference
\begin{equation}\label{eq:Hess-op}
\begin{split}
 A^u=A_0^u-a^u\colon\ H_1&\to H_0,
 \qquad
 \eta \mapsto B_u^{-1}\dot \eta+\bar L_u(\eta,\dot u) -a^u\eta
\end{split}
\end{equation}
of the para-Darboux Hessian operator $A_0^u$ in~(\ref{eq:Hess-op-0})
and the perturbation term $a^u$ in~(\ref{eq:a_u}).
For any $u\in U_2$, see~(\ref{eq:U_2}),
this operator is also bounded as a linear map
$
 A^u_2:=A^u|_{H_2}\colon H_2\to H_1
$.
\end{Definition}

\begin{Lemma}[$H_0$-symmetry]\label{le:a^u-symm}
At any $u\in U_1$, it holds
$\INNER{\xi}{A^u\eta}_{H_0}=\INNER{A^u\xi}{\eta}_{H_0}$
for all $\xi,\eta\in H_1$.
\end{Lemma}

\begin{proof}
Lemma~\ref{le:A^u_0-symm} and symmetry of the matrix~(\ref{eq:matrix-a_t}).
\end{proof}

\subsection{Fredholm operators}

\begin{Hypothesis}[on $u$]\label{hyp:u}
Let $u_-$ and $u_+$ be non-degenerate critical points of the
perturbed action $\Aa_h$.
Pick a \textit{basic path} $\hat u$ from $u_-$ to $u_+$
(see~\cite{Frauenfelder:2025e}),
i.e., $\hat u\in C^2(\R,U_2)$ with the property that there exists $T>0$
such that $\hat u(s)=u_-$ whenever $s\le -T$ and
$\hat u(s)=u_+$ whenever $s\ge T$.
Let $(H_0,H_1,H_2)$ be given by~(\ref{eq:triple}). Abbreviate
\begin{equation}\label{eq:abbreviate}
 W^{1,2}_{H_j}:=W^{1,2}(\R,H_j)
 ,\qquad
 L^2_{H_k}:=L^2(\R,H_k) .
\end{equation}
Let
\[
 u\in C^0(\R,U_1)
 ,\qquad
 u-\hat u\in W^{1,2}_{H_1}\cap L^2_{H_2} .
\]
\end{Hypothesis}

It can be shown that the non-degeneracy condition can be rephrased with the
help of the Hamiltonian flow $\varphi_t^h$, i.e., the flow of the
Hamiltonian vector field,
characterized by the requirement $\varphi_0^h=\id$ and
$\tfrac{{\rm d}}{{\rm d}t}\varphi_t^h=X_{h_t}\circ \varphi_t^h$.
Namely, $\Aa_h$ is Morse if and only if for every critical point $u$ of $\Aa_h$
we have $\ker \bigl({\rm d}\varphi_1^h u(0)-\1\bigr)=\{0\}$.

\begin{Theorem}\label{thm:para-Darboux}
For $u$ as in Hypothesis~$\ref{hyp:u}$
and $A^u$ as in~\eqref{eq:Hess-op}, the operators
\begin{equation*}
\begin{split}
 &\D^u
 =\p_s+A^u\colon\
 W^{1,2}_{H_0}\cap L^2_{H_1}
 \to L^2_{H_0},
\\
 &\D^u_2
 =\p_s+A^u_2\colon\
 W^{1,2}_{H_1}\cap L^2_{H_2}
 \to L^2_{H_1}
\end{split}
\end{equation*}
are both Fredholm operators of the same Fredholm index
\[
 \INDEX \D^u=\INDEX \D^u_2.
\]
\end{Theorem}

Theorem~\ref{thm:para-Darboux}
holds true by the abstract Theorem~\ref{thm:main}
which applies since the perturbed para-Darboux Hessian field $A$
is almost extendable by Theorem~\ref{thm:paradarboux-Hess-pert}.

\begin{Remark}[Conley--Zehnder index]
In the case where one of the asymptotic loops $u_-$ or $u_+$
is contractible, one can choose a filling disk for this loop.
This filling disk, together with the cylinder $u$,
induces a filling disk for the other asymptotic loop.
To each of these filling disks, one can associate a Conley--Zehnder
index by conjugating the linearized flow with a symplectic
trivialization over the filling disk.
In this case, one should be able to show
that the Fredholm index corresponds to the difference
of the Conley--Zehnder indices.
\end{Remark}

\begin{Remark}[linearized downward $L^2$ gradient flow]
Since $A^u$ is the $L^2$ Hessian of the functional $\Aa_h$
at $u\in U_1$ with respect to the standard inner product on $\R^{2n}$,
the kernel of the operator~$\D^u$ corresponds to the linearized downward
$L^2$ gradient flow of $\Aa_h$.
\end{Remark}

\subsubsection*{Main difficulties and how to overcome them}

Before introducing the abstract setup in Section~\ref{sec:alm-ext}
below which we use to prove
Theorem~\ref{thm:para-Darboux}, we first explain the main difficulties
and the main ideas how to overcome these difficulties.

The fact that $\D^u$ is Fredholm almost directly follows from
Rabier's theorem which itself generalizes previous work by
Robbin and Salamon.
In particular, Rabier does not need any differentiability assumption
on $s\mapsto A(s)$, but requires continuity.
However, since $u$ as a map $\R\to U_2$, $s\mapsto u_s$,
is only of class $L^2_{\rm loc}$, it is not necessarily continuous.
In particular, the map $\R\to \Ll(H_2,H_1)$, $s\mapsto A^{u_s}_2$,
is not necessarily of class $C^0$.
Thus the improvement of the Robbin--Salamon Fredholm
theorem~\cite{robbin:1995a} from $C^1$ to $C^0$
by Rabier~\cite{Rabier:2004a} cannot be applied to $\D^u_2$.

\textbf{Decomposition.}
In order to deal with this difficulty, we decompose the operators
$\D^u$ and~$\D^u_2$ into two parts.
To this end, we introduce the notation
\[
 A^{u_s}=
 \underbrace{B_{u_s}^{-1}\p_t-a^{u_s}}_{=:F^{u_s}}
 +\underbrace{\bar L_{u_s}(\cdot,\dot u_s)}_{=:C^{\dot u_s}}
 =F^{u_s}+C^{\dot u_s}
\]
and
\begin{equation*}
 \D^u\xi
 =\underbrace{(\p_s+F^u)}_{=:\F^u}\xi
 +\underbrace{\bigl[s\mapsto C^{\dot u_s}\xi_s\bigr]}_{=:
 M_{C^{\dot u}}\xi}
 =\F^u\xi+M_{C^{\dot u}}\xi .
\end{equation*}
The advantage of this decomposition is the following.
The map
\[
 \R\to\Ll(H_2,H_1),\qquad s\mapsto
 F^{u_s}_2=B_{u_s}^{-1}\p_t |_{H_2}-a^{u_s}_2
\]
is continuous since it avoids $\dot u_s$, while the non-continuous map
\[
 \R\to\Ll(H_2,H_1)
 ,\qquad
 s\mapsto C^{\dot u_s}_2=\bar L_{u_s}(\cdot,\dot u_s) |_{H_2}
\]
is of lower order since it avoids $\p_t$.

Since $C^{\dot u}$ is of lower order, we show, maybe after an
asymptotic correction, with the help of
Theorem~\ref{thm:perturbation}
that the multiplication operator $M_{C^{\dot u}}$ as a map
\smash{$W^{1,2}_{H_1}\cap L^2_{H_2}\to L^2_{H_1}$} is compact.
But Fredholm property and index are stable under compact perturbation.
Thus to show that
\begin{equation*}
\begin{split}
 \D^u_2\colon\
 W^{1,2}_{H_1}\cap L^2_{H_2}
 &\to L^2_{H_1},
\qquad
 \xi \mapsto \p_s\xi+A^u_2\xi
\end{split}
\end{equation*}
is Fredholm of the same Fredholm index as $\D^u$
is equivalent to showing this for
\begin{equation*}
\begin{split}
 \F^u_2=\p_s+F^u\colon\
 W^{1,2}_{H_1}\cap L^2_{H_2}
 &\to L^2_{H_1},
\qquad
 \xi \mapsto \p_s\xi+\underbrace{B_u^{-1}\p_t\xi-a^u\xi}_{=:F^u\xi} .
\end{split}
\end{equation*}

The fact that $\R\to\Ll(H_2,H_1)$, $s\mapsto F^{u_s}$, is continuous,
while $s\mapsto A^{u_s}$ is not, makes it easier to prove
the Fredholm property for $\F^u_2$ (Rabier's theorem applies)
than for $\D^u_2$ \big(it does not apply by non-continuity of
$s\mapsto C^{\dot u_s}$\big).
The disadvantage of $F^{u_s}$ when compared to $A^{u_s}$ is that, while~$A^{u_s}$ is $H_0$-symmetric for every $s$,
necessarily not so is $F^{u_s}$.
However, Rabier's theorem can as well deal with non-symmetric
operators.

\section{Almost extendability}\label{sec:alm-ext}

Consider a \textit{Hilbert space pair} $(H_0,H_1)$,
that is $H_0$ and $H_1$ are both infinite-dimensional Hilbert spaces
such that as sets $H_1\subset H_0$ and inclusion $\iota\colon H_1\INTO H_0$ is
compact and dense.

Then, as we explain in~\cite[Section~2]{Frauenfelder:2024c},
there exists an unbounded monotone function
\[
 \gf{h}=\gf{h}(H_0,H_1) \colon\ \mathbb{N} \to (0,\infty),
\]
called \textit{pair growth function}, such that the pair $(H_0,H_1)$
is isometric to the pair $\bigl(\ell^2,\ell^2_{\gf{h}}\bigr)$,
see~\cite[Appendix~A]{Frauenfelder:2024c}, and from now on we identify
the pairs
\[
 (H_0,H_1)=\bigl(\ell^2, \ell^2_\gf{h}\bigr)
 ,\qquad
 \gf{h}=\gf{h}(H_0,H_1) .
\]
Here $\ell^2_\gf{h}$ is defined as follows. In general,
for $\gf{f}\colon\mathbb{N}\to(0,\infty)$ unbounded monotone
\smash{$\ell^2_\gf{f}$} is the space of all sequences $x=(x_\nu)_{\nu \in \mathbb{N}}$ with
\smash{$\sum_{\nu=1}^\infty \gf{f}(\nu)x_\nu^2 <\infty$}.
The space \smash{$\ell^2_\gf{f}$} becomes a Hilbert space if we endow it
with the inner product
\begin{equation*}
 \langle x,y \rangle_\gf{f}=\sum_{\nu=1}^\infty \gf{f}(\nu)x_\nu y_\nu
 ,\qquad
 \Norm{x}_\gf{f}:=\sqrt{\langle x,x \rangle_\gf{f}} .
\end{equation*}
Therefore, we can define for every real number $r$ a Hilbert space
$ H_r:=\ell^2_{\gf{f}^r}$.
For each $s<r$, the inclusion $H_r\subset H_s$
is compact and dense.

The resulting triple of Hilbert spaces $H_0$, $H_1$, and $H_{r=2}$
is called a \textit{Hilbert space triple}, notation $(H_0,H_1,H_2)$.
Pick a Hilbert space triple $(H_0,H_1,H_2)$.

\subsection{Weak Hessian fields and decompositions}

Fix a Hilbert space triple $(H_0,H_1,H_2)$.
Let $U_1\subset H_1$ be open. Set $U_2:=U_1\cap H_2$.

\begin{Definition}[weak Hessian~\cite{Frauenfelder:2024c}]
An element $A$ of the Banach space\footnote{
 $\Ll(H_1,H_0)\cap \Ll(H_2,H_1)$ is a Banach space under the norm
 $\max\bigl\{\norm{\cdot}_{\Ll(H_1,H_0)},\norm{\cdot}_{\Ll(H_2,H_1)}\bigr\}$.
 }
$\Ll(H_1,H_0)\cap \Ll(H_2,H_1)$ is said a \textit{weak Hessian}
if it satisfies these two~axioms:
\begin{labeling}{\texttt{(symmetry)}}\itemsep=0pt
\item[\texttt{(symmetry)}]
 $\forall x,y\in H_1$: $\INNER{A x}{y}_0=\INNER{x}{A y}_0$
 called \textit{$H_0$-symmetry}.
\item[\texttt{(fredholm)}]
 $A\colon H_1\to H_0$,
 $A_2:=A|_{H_2}\colon H_2\to H_1$
 are Fredholm of index zero.
\end{labeling}
\end{Definition}

We are interested not in a single weak Hessian, but in fields of weak Hessians.

\begin{Definition}[weak Hessian field]\label{def:weak-Hessian-field}
A \textit{weak Hessian field} on $U_1$ is a continuous map
$A\in C^0(U_1,\Ll(H_1,H_0)) \cap C^0(U_2,\Ll(H_2,H_1))$,
notation $u\mapsto A^u$, satisfying the two conditions:
\begin{labeling}{\texttt{(Fredholm)}}\itemsep=0pt
\item[\texttt{(Symmetry)}]
 At any point $u\in U_1$ there is $H_0$-symmetry in the sense that
 \begin{equation}\label{eq:u-0-symmetric}
 \forall x,y\in H_1\colon \
 \INNER{A^u x}{y}_0=\INNER{x}{A^u y}_0 .
 \end{equation}
\item[\texttt{(Fredholm)}]
 $\forall u\in U_1$:
 $A^u\colon H_1\to H_0$ is Fredholm of index zero.

 $\forall u\in U_2$:
 $A^u_2\colon H_2\to H_1$ is Fredholm of index zero.
\end{labeling}
\end{Definition}

\begin{Remark}
Given a weak Hessian field $A$, then along level two
each operator $A^u$ is a weak Hessian, in symbols for each $u\in U_2$,
since $U_2\subset U_1$, the operator $A^u$ lies in
$\Ll(H_1,H_0)\cap \Ll(H_2,H_1)$ and satisfies \texttt{(symmetry)}
and \texttt{(fredholm)}.
\end{Remark}

\begin{Definition}[extendability]
We say that a weak Hessian field $A$ on $U_1$ \textit{extends}
if~$A$ extends to a continuous map
\smash{$U_1\to \Ll(H_1,H_0)\cap\underline{\Ll(H_2,H_1)}$},
still denoted by $u\mapsto A^u$, such that the restriction
$A^u_2:=A^u|_{H_2}\colon H_2\to H_1$ is Fredholm of index zero at
every point $u$ of \smash{\underline{$U_1$}}, and not only of $U_2$.
\end{Definition}

\begin{Remark}[equivalent formulation]
A weak Hessian field $A$ on $U_1$ \textit{extends}
if and only if~$A$ extends to a~continuous map
\smash{$U_1\to \Ll(H_1,H_0)\cap\underline{\Ll(H_2,H_1)}$}, still denoted by
$u\mapsto A^u$, such that~$A^u$
is a~weak Hessian at every point $u\in U_1$.
\end{Remark}

In general, the extendability condition is too strong as the
example of the area functional in a non-Darboux chart
shows, see~(\ref{eq:Hess-op-0}).\footnote{For $u\in U_1\subset W^{1,2}$ $\bigl(\dot u\in L^2\bigr)$ and
 $\eta\in W^{2,2}$, the product $\eta \dot u$ does not necessarily lie
 in $W^{1,2}$, however it does if $\dot u\in W^{1,2}$ $(u\in U_2)$;
 see~(\ref{eq:Behz-Hol-1}).}

\begin{Definition}[almost extendability]\label{def:alm.ext}\quad
\begin{enumerate}\itemsep=0pt
\item[(i)] We say that a weak Hessian field $A$ on $U_1$ \textit{almost extends}
if there exists a decomposition
\[
 A=F+C
\]
with
\begin{equation}\label{eq:hyp-F}
 F\in C^0(U_1,\Ll(H_1,H_0)\cap\Ll(H_2,H_1))
\end{equation}
and
\begin{equation}\label{eq:hyp-C}
 \exists r\in[0,1)\colon\
 C\in C^0\bigl(U_1,\Ll\bigl(\underline{H_r},H_0\bigr)\bigr)\cap C^0\bigl(U_2,\Ll\bigl(\underline{H_1},H_1\bigr)\bigr)
\end{equation}
such that the following two axioms hold:
\begin{labeling}{\texttt{(F)}}
\item[\texttt{(F)}]
 $\forall u \in\underline{\, U_1}\colon$
 $F^u_2:=F^u|_{H_2}\colon H_2\to H_1$
 is Fredholm of index zero.

\item[\texttt{(C)}]
 $\forall u\in U_1$ there exists an $H_1$-open neighborhood $V_u$ of $u$
 and a constant $\kappa$ such that
 for all $v,w\in V_u\cap H_2$ it holds the
 \emph{scale Lipschitz estimate}\footnote{Axiom~(\ref{eq:(C)}) is motivated by calculation~(\ref{eq:C-paraDarboux})
 for the para-Darboux Hessian.}
 \begin{equation}\itemsep=0pt\label{eq:(C)}
 \norm{C^v-C^w}_{\Ll(H_1)}
 \le\kappa\bigl(\abs{v-w}_{H_2}
 +\min\{\abs{v}_{H_2},\abs{w}_{H_2}\}\cdot\abs{v-w}_{H_1}\bigr).
 \end{equation}
\end{labeling}
\item[(ii)] If $A$ almost extends we call the pair $(F,C)$
a \textit{decomposition of $\boldsymbol{A}$}.
\end{enumerate}
\end{Definition}

\begin{Remark}[extendable $\Rightarrow$ almost extendable]
If a weak Hessian field $A$ is extendable,
then it is almost extendable (choose $F=A$ and $C=0$).
\end{Remark}

\begin{Lemma}\label{le:Eu-Fred}
Assume a weak Hessian field $A$ along $U_1$ is almost extendable
with decomposition~$(F,C)$. Then, for any $u\in U_1$,
the level one map $F^u\colon H_1\to H_0$ is Fredholm of
index zero, too.
\end{Lemma}

\begin{proof}
We write $F^u=A^u-C^u|_{H_1}$, where $A^u\colon H_1\to H_0$ is Fredholm of
index zero by axiom \texttt{(Fredholm)}.
Since $C^u\in\Ll(H_r,H_0)$, we can write $C^u|_{H_1}$
as a composition $C^u\circ\iota_r\colon {H_1\to H_r\to H_0}$
of a bounded and a compact operator. Hence $C^u|_{H_1}$ is compact.
By the stability of the Fredholm property and the index
under compact perturbation, we have that
$F^u$ is a Fredholm operator and
$\INDEX F^u=\INDEX A^u=0$.
This proves Lemma~\ref{le:Eu-Fred}.
\end{proof}

\subsection{Main theorem}

Let $(H_0,H_1,H_2)$ be a Hilbert space triple and $U_1\subset H_1$ an
open subset. Write
\[
 W^{1,2}_{H_j}:=W^{1,2}(\R,H_j)
 ,\qquad
 L^2_{H_j}:=L^2(\R,H_j) .
\]
To avoid constants,
maybe after replacing the norms by equivalent norms,
we assume in the following that
\[
 \abs{\cdot}_{H_0}\le \abs{\cdot}_{H_1}\le \abs{\cdot}_{H_2}.
\]

\begin{Definition}[connecting paths]
Fix two points $u_-,u_+\in U_2:=U_1\cap H_2$.
Fix a \textit{basic path} $\hat u$ from $u_-$ to $u_+$
(see~\cite{Frauenfelder:2025e}),
i.e., $\hat u\in C^2(\R,U_2)$ with the property that there exists $T>0$
such that $\hat u(s)=u_-$ whenever $s\le -T$ and
$\hat u(s)=u_+$ whenever $s\ge T$.

A \textit{connecting path from \boldmath$u_-$ to $u_+$}
is a continuous map $u\colon \R\to U_1$
such that the difference $u-\hat u$ lies in the intersection
Hilbert space \smash{$W^{1,2}_{H_1}\cap L^2_{H_2}$}, i.e.,
\begin{equation}\label{eq:connecting-paths}
 u\in C^0(\R,U_1)
 ,\qquad
 u-\hat u\in W^{1,2}_{H_1}\cap L^2_{H_2} .
\end{equation}
\end{Definition}

\begin{Remark}[independence of choice of basic path]
The notion of connecting path does not depend on the choice of
the basic path $\hat u$.
Indeed, suppose $\hat v$ is another basic path, then
$ u-\hat v
 =u-\hat u+\hat u-\hat v$
where \smash{$u-\hat u\in W^{1,2}_{H_1}\cap L^2_{H_2}$} by assumption
and the other difference is $C^2$ and of compact support, namely
\smash{$\hat u-\hat v\in C^2_{\rm c}(\R,H_2)\subset W^{1,2}_{H_1}\cap L^2_{H_2}$}.
\end{Remark}

\begin{Theorem}\label{thm:main}
Let $A$ be an almost extendable weak Hessian field on~$U_1$.
Consider two points $u_-,u_+\in U_2$ and a connecting path $u$.
Assume that both asymptotic operators
$A^{u_\mp}$ are isomorphisms as maps $H_1\to H_0$.
Then the operators
\begin{equation*}
 \D^u=\p_s+A^u\colon\
 W^{1,2}_{H_0}\cap L^2_{H_1}
 \to L^2_{H_0}
\end{equation*}
and
\begin{equation*}
 \D^u_2=\p_s+A^u_2\colon\
 W^{1,2}_{H_1}\cap L^2_{H_2}
 \to L^2_{H_1}
\end{equation*}
are both Fredholm operators of the same Fredholm index,
$\INDEX \D^u=\INDEX \D^u_2$.
\end{Theorem}

\begin{Remark}
That $A$ is almost extendable is only used to show that $\D^u_2$ is a
Fredholm operator.
That $\D^u$ is a Fredholm operator follows for every weak Hessian field
on $U_1$ satisfying the asymptotic non-degeneracy condition at $u_\mp$.
\end{Remark}

\begin{Lemma}\label{le:decomp}
Suppose that $A$ is an almost extendable weak Hessian field.
Let $u_*\in U_2$.
Then there exists a decomposition $A=F+C$ satisfying $C^{u_*}=0$.
\end{Lemma}

\begin{proof}
Since $A$ is extendable, there exists a decomposition $A=\tilde
F+\tilde C$.
In particular, there is $r\in[0,1)$ such that
$\tilde C\in C^0(U_1,\Ll(H_r,H_0))\cap C^0(U_2,\Ll(H_1))$.
Since $u_*\in U_2\subset U_1$, it holds that
$
 \tilde C^{u_*}\in \Ll(H_r,H_0)\cap \Ll(H_1)
$.

As $U_1$ is open, there exists $\eps>0$ such that the
ball \smash{$B_\eps^{H_1}(u_*)=\{v\in H_1\mid \abs{v-u_*}_{H_1}<\eps\}$} is contained
in $U_1$.
Choose a smooth cut-off function $\beta\colon\R\to[0,1]$
with the property that $\beta\equiv 1$ on $(-\infty,0]$
and $\beta\equiv 0$ on $\bigl[\eps^2,\infty\bigr)$
and define the map
\[
 v\mapsto
 P^v:=
 \begin{cases}
 \beta\bigl(\abs{u_*-v}^2_{H_1}\bigr) \tilde C^{u_*},
 &\text{$v\in B_{\eps}^{H_1}(u_*)$},
 \\
 0, &\text{else},
 \end{cases}
\]
which is element of $C^\infty(U_1,\Ll(H_r,H_0)\cap \Ll(H_1))$.
The map defined by
\[
 C:=\tilde C-P
\]
lies in $C^0(U_1,\Ll(H_r,H_0))\cap C^0(U_2,\Ll(H_1))$
since $\tilde C$ does and so does $P$,
because inclusion $U_2\INTO U_1$ is continuous.
Furthermore, the map defined by
\[
 F:=\tilde F+P
\]
lies in $C^0(U_1,\Ll(H_1,H_0)\cap\Ll(H_2,H_1))$ since $\tilde F$ does
and so does $P$, because the inclusions ${H_1\INTO H_r}$ and $H_2\INTO
H_1$ are continuous.
Note that
\[
 C^{u_*}=\tilde C^{u_*}-P^{u_*}
 =\tilde C^{u_*}-\beta(0) \tilde C^{u_*}=0
 ,\qquad
 F+C=\tilde F+\tilde C=A .
\]

\begin{Claim*}
The pair $(F,C)$ is a decomposition of $A$.
\end{Claim*}

To see this, we need to verify the two axioms
\texttt{(F)} and \texttt{(C)}.

\texttt{(F)} Let $u\in U_1$, then
\begin{equation*}
 F^u_2 =F^u|_{H_2}
 =\tilde F^u|_{H_2} +P^u|_{H_2}
 =\tilde F^u_2+\beta\bigl(\abs{u_*-u}_{H_1}^2\bigr)\tilde C^{u_*}|_{H_2} .
\end{equation*}
Since \smash{$\bigl(\tilde F,\tilde C\bigr)$} is a decomposition by axiom \texttt{(F)}
for \smash{$\tilde F$}, we have that \smash{$\tilde F^u_2\colon H_2\to H_1$}
is Fredholm of index zero.
Since \smash{$\tilde C^{u_*}\in \Ll(H_1)$} and
the inclusion $H_2\INTO H_1$ is compact,
it follows that \smash{$\tilde C^{u_*}|_{H_2}\colon H_1\to H_1$} is compact.
This proves \texttt{(F)} for $F$ since Fredholm property and index are
preserved under compact perturbation.

\texttt{(C)}
To show the local Lipschitz condition for $C$, we show that
$P$ is globally Lipschitz, then the local Lipschitz condition for
\smash{$C=\tilde C-P$} follows from the local Lipschitz condition for $\tilde
C$.
To this end, by the triangle inequality,
first we see that the function
\[
 U_1\to\R
 ,\qquad
 v\mapsto\abs{u_*-v}_{H_1}
\]
is Lipschitz.\footnote{
 Add $0=-w+w$ to $v$, then
 $\abs{u_*-v}_{H_1}-\abs{u_*-w}_{H_1}
 \le \abs{u_*-w}_{H_1}+\abs{w-v}_{H_1}-\abs{u_*-w}_{H_1}$
 by the triangle inequality, interchange $v$ and $w$ to get
 $\abs{u_*-w}_{H_1}-\abs{u_*-v}_{H_1}\le \abs{v-w}_{H_1}$.
 }
Since the squaring function $\R\ni x\mapsto x^2$ on the compact set
$[0,\eps]$ is Lipschitz and~$\beta$ is Lipschitz, the function
\[
 U_1\to\R
 ,\qquad
 v\mapsto\beta\bigl(\abs{u_*-v}_{H_1}^2\bigr)
\]
is Lipschitz (Lipschitz is preserved under composition).
As $\tilde C^{u_*}$ does not depend on $v$, the map
\[
 U_1\to \Ll(H_r,H_0)\cap\Ll(H_1)
 ,\qquad
 v\mapsto P^v
\]
is Lipschitz, say with Lipschitz constant $L$. That is, for all
$v,w\in U_1$, it holds
\[
 \norm{P^v-P^w}_{\Ll(H_r,H_1)\cap\Ll(H_1)}\le L\abs{v-w}_{H_1} .
\]
If $v,w\in U_2$, this implies the estimate
\[
 \norm{P^v-P^w}_{\Ll(H_1)}
 \le \norm{P^v-P^w}_{\Ll(H_r,H_1)\cap\Ll(H_1)}
 \le L\abs{v-w}_{H_1}
 \le L\abs{v-w}_{H_2}.
\]
Hence $P\colon U_2\to\Ll(H_1)$ is Lipschitz continuous
even globally on $U_2$.
This proves axiom~\texttt{(C)} for~$C$, hence the claim.

The proof of Lemma~\ref{le:decomp} is complete.
\end{proof}

\begin{Remark}[nullifying $C$ on finitely many points]
\label{rem:decomp}
The proof of Lemma~\ref{le:decomp} shows that one can improve the
statement of the lemma to a stronger statement:
Assume $A$ is an almost extendable weak Hessian field on $U_1$
and $\Delta\subset U_2$ is a finite subset.
Then there is a decomposition $A=F+C$
such that $C$ vanishes along the points of $\Delta$.
The case $\Delta=\{u_-,u_+\}$ matters.
\end{Remark}

\subsection{Proof of main theorem}

The proof of the main theorem, Theorem~\ref{thm:main},
has three parts.

\begin{PartI}
The operator
\smash{$\D^u=\p_s+A^u\colon W^{1,2}_{H_0}\cap L^2_{H_1}\to L^2_{H_0}$}
is Fredholm.
\end{PartI}

\begin{proof}
Verifying (H1)--(H3) and applying Corollary~\ref{cor:symm-Fredholm} proves Part I.
\begin{enumerate}\itemsep=0pt
\item[(H1)] True since $(H_0,H_1)$ is a Hilbert space pair.
\item[(H2)] Since $W^{1,2}(\R,H_1)$ embeds in $C^0 (\R,H_1)$,
see, e.g.,~\cite[Appendix~A.6]{Frauenfelder:2025e},
the map $\R\ni s\mapsto u_s\in U_1\subset H_1$ is continuous.
Since moreover $A\colon U_1\to\Ll(H_1,H_0)$ is
continuous as well, the map
$\R\to \Ll(H_1,H_0)$, $s\mapsto A^{u_s}$
is continuous. Hence (H2) is satisfied.

\item[(H3)] By~(\ref{eq:connecting-paths}), the path $u$
converges to $u_\mp$, as $s\to\mp\infty$, hence
$A^{u_s}$ converges to operators $A^{u_\mp}$, as $s\to\mp\infty$,
and as maps $H_1\to H_0$ these are invertible by assumption.
Hence (H3) is satisfied. This proves Part~I.
\hfill\qed
\end{enumerate}
\renewcommand{\qed}{}
\end{proof}

\begin{PartII}The operator
$\D^u_2=\p_s+A^u_2\colon W^{1,2}_{H_1}\cap L^2_{H_2}\to L^2_{H_1}$
is Fredholm.
\end{PartII}

\begin{proof}
By assumption, $A^u$ is almost extendable,
so there exists a decomposition $A^u=F^u+C^u$.
By Remark~\ref{rem:decomp}, we can assume in addition
that $C^{u_-}=0=C^{u_+}$.

By invariance of the Fredholm property under compact perturbation
and since $\D^u_2=\p_s +F^u_2+C^u_2$, to prove Part~II it suffices to
show the following:
\begin{itemize}\setlength{\leftskip}{0.10cm}\itemsep=0pt
\item[\bf II-a.]
 \textit{The operator
 \smash{$\F^u_2=\p_s+F^u_2\colon W^{1,2}_{H_1}\cap L^2_{H_2}\to L^2_{H_1}$}
 is Fredholm.}
\item[\bf II-b.]
 \textit{The multiplication operator
 \smash{$M^{C^u}_2\colon W^{1,2}_{H_1}\cap L^2_{H_2}\to L^2_{H_1}$}
 is compact.}
\end{itemize}

\begin{proof}[Proof of II-a.]
It suffices to verify hypotheses (H1)--(H5)
in Rabier's Theorem~\ref{thm:Rabier}.
\begin{enumerate}\itemsep=0pt
\item[(H1)] The pair $(H,W):=(H_1,H_2)$ is a Hilbert space pair.

\item[(H2)] The map
$\R\to U_1 \to \Ll(H_2,H_1)$, $s\mapsto u_s\mapsto F^{u_s}_2$
is a composition of continuous maps, by~(\ref{eq:connecting-paths})
and~(\ref{eq:hyp-F}), hence continuous.

\item[(H3)]
Since $u$ is a connecting path from $u_-$ to $u_+$, we have
that $\lim_{s\to\mp\infty} u_s=u_\mp$.
Since $F\in C^0(U_1,\Ll(H_2,H_0))$, by~(\ref{eq:hyp-F}),
it holds that
$\lim_{s\to\mp\infty}\norm{F^{u_s}-F^{u_\mp}}_{\Ll(H_2,H_1)}
 =0$.
It remains to show that the asymptotic
limit operators $F^{u_\mp}$ are isomorphisms as linear maps $H_2\to H_1$.
To see this, observe that $F^{u_\mp}=A^{u_\mp}$, since $C^{u_\mp}=0$ as
mentioned above, and that by hypothesis $A^{u_\mp}\colon H_1\to H_0$
are isomorphisms.
In particular, the kernel vanishes, hence so does the kernel of the
restriction $\ker A^{u_\mp}_2\colon H_2\to H_1$.
Since $u_\mp$ are elements of $U_2$, by axiom \texttt{(Fredholm)}
of weak Hessian field, the operator $A^{u_\mp}_2\colon H_2\to H_1$ is
Fredholm of index zero.
Hence injective is equivalent to surjective
and therefore by the open mapping theorem the two maps
\smash{$A^{u_\mp}_2$} are isomorphisms. But \smash{$F^{u_\mp}_2=A^{u_\mp}_2$}.

\item[(H4)] For $F^{u_s}_2$ whenever $s\in\R$:
By Lemma~\ref{le:Rabier-symmetric}, based on
Definition~\ref{def:weak-Hessian-field} for~$A$,
we know that $A^{u_s}\colon H_1\to H_0$ satisfies~(H4).
Our first goal is to show that
$F^{u_s} =A^{u_s}-C^{u_s}|_{H_1}\colon H_1\to H_0$ satisfies (H4).
To see this, observe that $C^{u_s}|_{H_1}\colon H_1\to H_0$ is a compact
operator: indeed, $C^{u_s}\in\Ll(H_r,H_0)$ by~(\ref{eq:hyp-C})
and inclusion $H_r\INTO H_1$ is compact since $r<1$.
But (H4) is stable under compact perturbation
by Proposition~\ref{prop:Rabier+compact}.
Thus $F^{u_s}\colon H_1\to H_0$ satisfies (H4)
whenever $s\in\R$ proving the goal.
Now by Lemma~\ref{le:Rabier+level-2},
based on axiom~\texttt{(F)},
the restriction $F^{u_s}_2\colon H_2\to H_1$
satisfies (H4) as well.

\item[(H5)]
Since $F^{u_\mp}=A^{u_\mp}$, as we saw during the verification of
(H3), and the spectrum of $A^{u_\mp}$ is real but does not contain
zero, by invertibility, it holds that
${\rm i}\R\cap\spec A^{u_\mp}=\varnothing$.
\end{enumerate}

Rabier's Theorem~\ref{thm:Rabier} concludes the proof of Step~II-a.
\end{proof}

\begin{proof}[Proof of II-b.]
By Theorem~\ref{thm:perturbation},
it suffices to show that $C^u\in L^2(\R,\Ll(H_1))$.

\textit{Basic path.}
Let $\hat u$ be a basic path from $u_-$ to $u_+$
with interval, say $[-T, T]$.

\textit{Positive end.}
By axiom~(\texttt{C}), for $C$
there is an $H_1$-open neighborhood $V_{u_+}$ of $u_+$,
and a~constant $\kappa_+>0$ such that all $v,w\in V_{u_+}\cap H_2$
satisfy the scale Lipschitz estimate~(\ref{eq:(C)}).
In particular, since $w:=u_+\in U_2$ and $C^{u_+}=0$,
for any $v\in V_{u_+}\cap H_2$ there is the local Lipschitz estimate
\begin{align}
 \norm{C^v}_{\Ll(H_1)}
 &\le\kappa\Bigl(\abs{v-u_+}_{H_2}
 +\underbrace{\min\{\abs{v}_{H_2},\abs{u_+}_{H_2}\}}_{\le\abs{u_+}_{H_2}}
 \cdot \underbrace{\abs{v-u_+}_{H_1}}_{\le \abs{v-u_+}_{H_2}}\Bigr)
 \le\kappa_+ \abs{v-u_+}_{H_2},\label{eq:pos-end}
\end{align}
where $\kappa_+=\kappa_+(u_+)=\kappa+\abs{u_+}_{H_2}$.
Since the connecting path $s\mapsto u_s$ converges to $u_+$
in $H_1$, there exists a time $T_+\ge T$ such that
$u_s\in V_{u_+}$ whenever $s\ge T_+$.

\textit{Negative end.}
Similarly, there exist $V_{u_-}$, $\kappa_->0$, $T_-\ge T$ such that
\begin{equation}\label{eq:neg-end}
 \norm{C^v}_{\Ll(H_1)}
 \le\kappa_-\abs{v-u_-}_{H_2}
\end{equation}
whenever $v\in V_{u_-}\cap H_2$
and $u_s\in V_{u_-}$ whenever $s\le -T_-$.

\textit{Compact part.}
For every $s\in[-T_-,T_+]$, by axiom~(\texttt{C}),
we choose an $H_1$-open neighborhood $V_s$ of $u_s$ in $U_1$
such that there exists a constant $\kappa_s>0$ with the property
that for all $v,w\in V_s\cap H_2$ it holds the scale Lipschitz estimate
\begin{equation}\label{eq:hhjhj777bb}
 \norm{C^v-C^w}_{\Ll(H_1)}
 \le\kappa_s (\abs{v-w}_{H_2}
 +\min\{\abs{v}_{H_2},\abs{w}_{H_2}\}\cdot\abs{v-w}_{H_1} ).
\end{equation}
Since $u\colon\R\to U_1$ is continuous and the interval
$[-T_-,T_+]$ is compact,
the image $u_{[-T_-,T_+]}$ is a~compact subset of $U_1$.
Therefore, there exists a positive integer $N$ and times
$
 -T_-<s_1<s_2<\dots<s_N<T_+
$
such that the $V_{s_j}$ cover the image
\[
 u_{[-T_-,T_+]}\subset \bigcup_{j=1}^N V_{s_j} .
\]
As $H_2$ is dense in $H_1$, for any $j=1,\dots,N$,
we can choose $p_j\in V_{s_j}\cap H_2$
as illustrated by Figure~\ref{fig:fig-covering-C^u}.

\begin{figure}[!ht]
 \centering
 \includegraphics
 {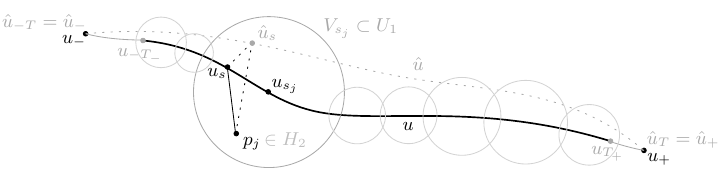}
 \caption{Basic path $\hat u$ and open cover of
 connecting path $u$ along time $[-T_-,T_+]$.}
 \label{fig:fig-covering-C^u}
\end{figure}

\textit{Estimates.}
Firstly, we assume that $s\in[-T_-,T_+]$.
Then $u_s$ lies in one of the finitely many
open sets, i.e., there exists $j=j(s)$ such that
$u_s\in V_{s_j}$. We~estimate
\begin{equation*}
\begin{split}
 \norm{C^{u_s}}_{\Ll(H_1)}
 &\le \norm{C^{p_j}}_{\Ll(H_1)}+\norm{C^{u_s}-C^{p_j}}_{\Ll(H_1)}
\\
 &\le \norm{C^{p_j}}_{\Ll(H_1)}
 +\kappa_{s_j}
 \Bigl(\abs{u_s-p_{j}}_{H_2}
 +\underbrace{\min\{\abs{u_s}_{H_2},\abs{p_{j}}_{H_2}\}}_{\le \abs{p_{j}}_{H_2}}
 \cdot \underbrace{\abs{u_s-p_{j}}_{H_1}}_{\le \abs{u_s-p_{j}}_{H_2}}\Bigr)
\\
 &\le\underbrace{\max_{\ell=1,\dots,N}\norm{C^{p_\ell}}_{\Ll(H_1)}}_{=:a_u}
 +\underbrace{\max_{\ell=1,\dots,N}
 (\kappa_{s_\ell}+\abs{p_{\ell}}_{H_2} )}_{=:b_u}
 (\abs{u_s-\hat u_s}_{H_2}+\abs{\hat u_s-p_{j}}_{H_2} )
\\
 &\le
 \underbrace{a_u+b_u \max_{\sigma\in[-T,T]\atop \ell=1,\dots,N}
 \abs{\hat u_\sigma-p_{\ell}}_{H_2}}_{=:c_u}
 +b_u \abs{u_s-\hat u_s}_{H_2} .
\end{split}
\end{equation*}
Inequality~1 is by adding zero and the triangle inequality.
Inequality~2 uses~(\ref{eq:hhjhj777bb}) for $s=s_j$.
Inequality~3 is by taking maxima and adding again zero
$\hat u_s-\hat u_s$, as illustrated by
Figure~\ref{fig:fig-covering-C^u}.
Inequality~4 uses that a basic path $\hat u$ is continuous as a map to
$U_2$ and constant outside $[-T,T]$, so the maximum exists.

Secondly, for $s\ge T_+$, by~(\ref{eq:pos-end}) it holds
$ \norm{C^{u_s}}_{\Ll(H_1)}
 \le\kappa_+\abs{u_s-u_+}_{H_2}$.

Thirdly, for $s\le -T_-$, by~(\ref{eq:neg-end}) it holds
$ \norm{C^{u_s}}_{\Ll(H_1)}
 \le\kappa_-\abs{u_s-u_-}_{H_2}$.

To prove that $C^u\in L^2(\R,\Ll(H_1))$, we estimate
\begin{equation*}
\begin{split}
 \norm{C^u}_{L^2(\R,\Ll(H_1))}^2
 &=\biggl(\int_{-\infty}^{-T_-} +\int_{-T_-}^{T_+} +\int_{T_+}^\infty\biggr)
 \norm{C^{u_s}}_{\Ll(H_1)}^2 \, {\rm d}s
\\
 &\le\bigl(\kappa_-^2 +\kappa_+^2\bigr)\norm{u-u_-}_{L^2_{H_2}}^2
 +(T_-+T_+) 2 c_u^2+b_u\norm{u-\hat u}_{L^2([-T,T],H_2)}^2
\end{split}
\end{equation*}
which is finite.
This proves Step~II-b.
\end{proof}

Since \smash{$\D^u_2=\F^u_2+M^{C^u}_2$} and since Fredholm
property and index are invariant under compact perturbation, we get
from Steps~II-a and II-b that $\D^u_2$ is a Fredholm operator and that
\begin{equation*}
 \INDEX \D^u_2=\INDEX \F^u_2 .
\end{equation*}
This concludes the proof of Part~II.
\end{proof}

\begin{PartIII} Equal Fredholm indices $\INDEX \D^u=\INDEX \D^u_2$.
\end{PartIII}

\begin{proof}[Proof of Part III -- using a homotopy argument]
Because the space $C^0_{\rm c}(\R,H_2)$ is dense in $L^2(\R,H_2)$,
see, e.g.,~\cite[Theorem~A.5.4]{Frauenfelder:2025e},
the connecting path $u\in C^0(\R,U_1)$ from $u_+$ to $u_-$
can be approximated by a connecting path
which additionally is
continuous as a map to $H_2$ and satisfies
that the limit as $s\to\mp\infty$ exists in $H_2$ and is $u_{\mp}$.
In~particular, after a homotopy inside the space of connecting paths from
$u_-$ to $u_+$, we can assume that our connecting path $u$ is
additionally continuous on level $H_2$.

Since the Fredholm index is homotopy invariant,
it suffices to show the equality of indices only for connecting paths
which are additionally continuous on level $H_2$.
For such a path, the map $s\mapsto A^{u_s}_2$ is continuous as a map
$\R\to\Ll(H_2,H_1)$.
Since both operators $A^{u_s}\colon H_1\to H_0$ and
$A^{u_s}_2\colon H_2\to H_1$ are Fredholm of index zero, \cite[Corollary~3.4]{Frauenfelder:2024c}
shows that the two operators have the same spectrum, in particular
they have the same spectral flow.
By~\cite[Theorem~A]{Frauenfelder:2024d},
the Fredholm index of $\D^u=\p_s+A^u$ is the spectral flow of
$s\mapsto A^{u_s}$
and
the Fredholm index of~${\D^u_2=\p_s+A^u_2}$ is the spectral flow of
$s\mapsto A^{u_s}_2$.
Since the spectral flows agree, the two indices agree.
This proves Part~III.
\end{proof}

The proof of Theorem~\ref{thm:main} is complete.

\section{Para-Darboux Hessian is almost extendable}

Recall the setup from Section~\ref{sec:unp-action},
where an open subset $\mathfrak{U}\subset\R^{2n}$
carries an exact symplectic form $\omega={\rm d}\lambda$.
The Hilbert space triple $(H_0,H_1,H_2)$ is given by the
$W^{k,2}\bigl(\SS^1,\R^{2n}\bigr)$-Sobolev spaces~(\ref{eq:triple}).
Open subsets $U_1\subset H_1$ and $U_2\subset H_2$
are defined by~(\ref{eq:U_2}).
Recall further that the identity~(\ref{eq:B})
determines a smooth map
\[
 B\colon\ \mathfrak{U}\to\GL(2n,\R)
 ,\qquad
 x\mapsto B_x .
\]
In~(\ref{eq:L-(1,2)}), we defined a (1,2) tensor
$\bar L\colon \Xx(\mathfrak{U}) \times \Xx(\mathfrak{U})
\to \Xx(\mathfrak{U})$ where $\Xx(\mathfrak{U})$ is the set of vector
fields along $\mathfrak{U}$.

In order to prove the lemma and the theorem below,
we need the following result from fractional Sobolev theory.

\begin{Theorem}[{\cite[Theorem~7.4]{Behzadan:2021a}}]
Assume $\rho_1$, $\rho_2$, $\rho$ are real numbers satisfying
\begin{enumerate}\itemsep=0pt
\item[$(1)$]
 $\rho_1\ge \rho\ge 0$ and $\rho_2\ge \rho\ge 0$;
\item[$(2)$]
 $\rho_1+\rho_2>\frac12 +\rho$.
\end{enumerate}
Then the following is true.
If \smash{$v\in W^{\rho_1,2}\bigl(\SS^1\bigr)$} and \smash{$w\in W^{\rho_2,2}\bigl(\SS^1\bigr)$}, then
\smash{$vw\in W^{\rho,2}\bigl(\SS^1\bigr)$} and pointwise multiplication of functions is a
continuous bi-linear map
\begin{equation*}
 W^{\rho_1,2}(\SS^1)\times W^{\rho_2,2}\bigl(\SS^1\bigr)\to W^{\rho,2}\bigl(\SS^1\bigr) .
\end{equation*}
\end{Theorem}

Most important for us is the theorem in the forms
\begin{equation}\label{eq:Behz-Hol-1}
 W^{1,2}\times L^2\to L^2
 ,\qquad
 W^{2,2}\times W^{1,2}\to W^{1,2}
 ,
\end{equation}
and
\begin{equation}\label{eq:Behz-Hol-2}
 W^{r,2}\times L^2\stackrel{r>\frac12}{\longrightarrow} L^2
 ,\qquad
 W^{1,2}\times W^{1,2}\to W^{1,2}
 .
\end{equation}

\subsection{Unperturbed case}

\begin{Lemma}\label{le:para-weak}
The para-Darboux Hessians, one for each $u\in U_1$, defined by
\[
 A_0^u\colon \ H_1\to H_0
 ,\qquad
 \eta\mapsto B_u^{-1}\dot \eta+\bar L_u(\eta,\dot u)
\]
determine a weak Hessian field $A_0$ on $U_1$.
\end{Lemma}

\begin{proof}
By~(\ref{eq:Behz-Hol-1}),
the map $u\mapsto A^u_0$ is
element of the space $C^0(U_1,\Ll(H_1,H_0))$ and
of the space $C^0(U_2,\Ll(H_2,H_1))$.
The (\texttt{Symmetry}) axiom~(\ref{eq:u-0-symmetric}) holds true by
Lemma~\ref{le:A^u_0-symm}.
It remains to check the (\texttt{Fredholm}) axiom.
By Theorem~\ref{thm:Fredholm-F_C}
for the loop $t\mapsto \omega_{u_t}$, the operator as a map
\begin{equation}\label{eq:hjghjghj77799}
 B_u^{-1}\p_t\colon\ H_1\to H_0,\qquad H_2\to H_1
 ,\qquad \forall u\in U_1
\end{equation}
is Fredholm of index zero.
Using the Sobolev estimate $\abs{\eta}_{L^\infty}\le\abs{\eta}_{H_1}$,
the term
\[
 C^u\colon\ \eta\mapsto \bigl[t\mapsto \bar L_{u_t}(\eta_t,\dot u_t) \bigr]
\]
is bounded as a map $H_0\to H_0$ if $u\in U_1$ and as a map
$H_1\to H_1$ if $u\in U_2$.
So, by compactness of the embeddings
$H_1\INTO H_0$ and $H_2\INTO H_1$, both operators $C^u$ are compact.
But Fredholm property and index are stable under compact perturbation.
This proves Lemma~\ref{le:para-weak}.
\end{proof}

\begin{Theorem}\label{thm:paradarboux-Hess}
The para-Darboux weak Hessian field $A_0$ is almost extendable.
Moreover, the pair $(F,C)$ defined for $u\in U_1$ by
\[
 F^u\colon\ \eta\mapsto B_u^{-1}\dot \eta
 ,\qquad
 C^u\colon\ \eta\mapsto \bar L_u(\eta,\dot u) ,
\]
is a decomposition.
\end{Theorem}

\begin{proof}[Proof of Theorem~\ref{thm:paradarboux-Hess}]
The proof has four Steps~1, 2, (\texttt{C}), and (\texttt{F}).

\medskip\noindent
\textbf{Step~1.}
By~(\ref{eq:Behz-Hol-2}),
the map $u\mapsto C^u$ is element of the space
$C^0(U_1,\Ll(\underline{H_r},H_0))$ and of the space
$C^0(U_2,\Ll(\underline{H_1},H_1))$,
where $H_r=W^{r,2}\bigl(\SS^1,\R^{2n}\bigr)$.

\medskip\noindent
\textbf{Step~2.}
We need to show that the map $F\colon u\mapsto F^u=B_u^{-1}\p_t$
is element of the space $C^0(U_1,\Ll(H_1,H_0)\cap\Ll(H_2,H_1))$;
see~(\ref{eq:hyp-F}).

\begin{proof}
We prove that $F\in C^0(U_1,\Ll(H_2,H_1))$;
similarly $F\in C^0(U_1,\Ll(H_1,H_0))$.

Fix $u\in U_1$ and let $K=K(u)$ be a compact neighborhood
of the (compact) image of $u$ in~${\mathfrak{U}\subset\R^{2n}}$.
Pick $v\in U_1$ near $u$ taking values in $K$ as well.
The smooth map $x\mapsto B_x^{-1}$ is Lipschitz continuous
since $K=K(u)$ is compact. Let $\lambda_u$ be the Lipschitz constant.
Now, given any~${\xi\in H_2}$, we estimate
\begin{equation*}
\begin{split}
 \abs{(F^u-F^v)\xi}_{W^{1,2}}^2
 &=\bigl|\bigl(B_u^{-1}-B_v^{-1}\bigr)\dot\xi\bigr|_{L^2}^2
 +\bigl|
 \bigl(B_u^{-1}-B_v^{-1}\bigr)\ddot\xi
 +\bigl(\p_t\bigl(B_u^{-1}-B_v^{-1}\bigr)\bigr)\dot\xi
 \bigr|_{L^2}^2
\\
 &\le\lambda_u^2\abs{u-v}_{L^\infty}^2\abs{\xi}_{W^{1,2}}^2
 +2\lambda_u^2\abs{u-v}_{L^\infty}^2\abs{\xi}_{W^{2,2}}^2
 \\
 &\quad{}
 +4 \bigl|{\rm d}B_u^{-1}\bigl(\dot u\,\underline{-\,\dot v},\dot\xi\bigr)\bigr|_{L^2}^2
 +4 \bigl|\bigl(\underline{{\rm d}B_u^{-1}}-{\rm d}B_v^{-1}\bigr)\bigl(\dot v,\dot\xi\bigr)\bigr|_{L^2}^2
\\
 &\le3\lambda_u^2\abs{u-v}_{W^{1,2}}^2\abs{\xi}_{W^{2,2}}^2
 +4\bigl|{\rm d}B_u^{-1}\bigr|_{L^\infty(K)}^2\abs{u-v}_{W^{1,2}}^2\bigl|\dot\xi\bigr|_{L^\infty}^2
 \\
 &\quad{}
 +4\lambda_u^2\abs{u-v}_{L^\infty}^2\abs{v}_{W^{1,2}}^2\bigl|\dot\xi\bigr|_{L^\infty}^2
\\
 &\le\bigl(3\lambda_u^2+4\bigl|{\rm d}B_u^{-1}\bigr|_{L^\infty(K)}^2
 +4\lambda_u^2 \abs{v}_{W^{1,2}}^2\bigr)
 \abs{u-v}_{W^{1,2}}^2 \abs{\xi}_{W^{2,2}}^2.
\end{split}
\end{equation*}
Here inequality one uses Lipschitz continuity of $K\ni x\mapsto B_x^{-1}$
with Lipschitz constant $\lambda_u$. We also added \textit{zero}
and used the triangle inequality.
Inequality two uses the Sobolev embedding $W^{1,2}\INTO C^0$ with
constant $1$, we also used Lipschitz continuity again.
This proves Step~2.
\end{proof}

\noindent
\textbf{Step~(\texttt{C}).} The map $u\mapsto C^u$ satisfies the scale
Lipschitz estimate~(\ref{eq:(C)}).

\begin{proof}
It remains to check the local scale Lipschitz axiom~(\texttt{C}),
see~(\ref{eq:(C)}), for the map~$C$.
To this end suppose $u\in U_1=W^{1,2}\bigl(\SS^1,\mathfrak{U}\bigr)$.
Since $W^{1,2}\subset C^0$, the image of $u$ is a compact subset
of~$\mathfrak{U}$. Therefore, there exist, firstly, an open subset
$\mathfrak{V}$ which contains the image of $u$ and, secondly, a~constant $c$ such that
\begin{align}
 &\bigl\|\bar L_x\bigr\|_{\Ll(\R^{2n}\times\R^{2n};\R^{2n})}
\le c, \nonumber
\\
 &\bigl\|\bar L_x-\bar L_x\bigr\|_{\Ll(\R^{2n}\times\R^{2n};\R^{2n})}
\le c \abs{x-y}_{\R^{2n}}, \nonumber
\\
 &\bigl\| {\rm d}\bar L_x\bigr\|_{\Ll(\R^{2n}\times\R^{2n}\times\R^{2n};\R^{2n})}
\le c ,\nonumber
\\
 &\bigl\|{\rm d}\bar L_x-{\rm d}\bar L_x\bigr\|_{\Ll(\R^{2n}\times\R^{2n}\times\R^{2n};\R^{2n})}
\le c \abs{x-y}_{\R^{2n}},\label{eq:L-est}
\end{align}
whenever $x,y\in\mathfrak{V}$.
We define an open neighborhood of $u$ in $U_1$ by
\[
 V_u:=W^{1,2}\bigl(\SS^1,\mathfrak{V}\bigr)
 \cap B_1^{H_1}(u),
\]
where \smash{$B_1^{H_1}(u)$} is the open radius $1$ ball in $H_1$ centered at $u$.
Pick elements $v,w\in V_u\cap H_2 =V_u\cap W^{2,2}$.
Note that
\[
 \abs{v}_{W^{1,2}}\le \abs{u}_{W^{1,2}}+1
 ,\qquad
 \abs{w}_{W^{1,2}}\le \abs{u}_{W^{1,2}}+1 .
\]
Now for $\xi\in H_1$, we estimate
\begin{equation*}
\begin{split}
 \abs{(C^v-C^w)\xi}_{W^{1,2}}^2
 \le{}& \bigl|\bar L_v(\xi,\dot v)-\bar L_w(\xi,\dot w)\bigr|_{L^2}^2
 +\bigl|\p_t\bigl(\bar L_v(\xi,\dot v)-\bar L_w(\xi,\dot w)\bigr)\bigr|_{L^2}^2
\\
 \le{}& 2\bigl|\bar L_v(\xi,\dot v-\dot w)\bigr|_{L^2}^2
 +2\bigl|\bigl(\bar L_v-\bar L_w\bigr) (\xi,\dot w)\bigr|_{L^2}^2
 \\
 &{+}\,7\bigl|{\rm d}\bar L_v(\dot v-\dot w,\xi,\dot v)\bigr|_{L^2}^2
 +7\bigl|{\rm d}\bar L_v(\dot w,\xi,\dot v-\dot w)\bigr|_{L^2}^2
 \\
 &{+}\,7\bigl|\bigl({\rm d}\bar L_v-{\rm d}\bar L_w\bigr)(\dot w,\xi,\dot w)\bigr|_{L^2}^2
 \\
 &{+}\,7\bigl|\bar L_v\bigl(\dot \xi,\dot v-\dot w\bigr)\bigr|_{L^2}^2
 +7\bigl|\bigl(\bar L_v-\bar L_w\bigr)\bigl(\dot \xi,\dot w\bigr)\bigr|_{L^2}^2
 \\
 &{+}\,7\bigl|\bar L_v(\xi,\ddot v-\ddot w)\bigr|_{L^2}^2
 +7\bigl|\bigl(\bar L_v-\bar L_w\bigr)(\xi,\ddot w)\bigr|_{L^2}^2,
\end{split}
\end{equation*}
where in inequality two we added four times zero
and used Lemma~\ref{le:sum-quad} for $k=7$ summands.
Now we estimate summand by summand starting below
\begin{equation*}
\begin{split}
 \bigl|\bigl(\bar L_v-\bar L_w\bigr)(\xi,\ddot w)\bigr|_{L^2}^2
 &=\int_0^1\bigl|\bigl(\bar L_{v_t}-\bar L_{w_t}\bigr) (\xi_t,\ddot w_t)\bigr|^2\, {\rm d}t
\\
 &\le \int_0^1c^2
 \abs{v_t-w_t}^2\abs{\xi_t}^2\abs{\ddot w_t}^2\, {\rm d}t
\\
 &\le c^2 \abs{v-w}_{L^\infty}^2\abs{\xi}_{L^\infty}^2
 \abs{\ddot w}_{L^2}^2
\\
 &\le c^2 \abs{v-w}_{W^{1,2}}^2\abs{\xi}_{W^{1,2}}^2
 \abs{w}_{W^{2,2}}^2 .
\end{split}
\end{equation*}
Here inequality one uses~(\ref{eq:L-est}) and the Sobolev embedding
$L^\infty\INTO W^{1,2}$ with constant $1$.
Similarly, we estimate
\begin{equation*}
\begin{split}
& \bigl|\bar L_v(\xi,\ddot v-\ddot w)\bigr|_{L^2}^2
\le c^2 \abs{\xi}_{W^{1,2}}^2\abs{v-w}_{W^{2,2}}^2,
\\
& \bigl|\bigl(\bar L_v-\bar L_w\bigr)\bigl(\dot \xi,\dot w\bigr)\bigr|_{L^2}^2
\le c^2 \abs{v-w}_{W^{1,2}}^2 \abs{\xi}_{W^{1,2}}^2 \abs{w}_{W^{2,2}}^2,
\\
& \bigl|\bar L_v\bigl(\dot \xi,\dot v-\dot w\bigr)\bigr|_{L^2}^2
\le c^2 \abs{\xi}_{W^{1,2}}^2 \abs{v-w}_{W^{2,2}}^2,
\\
& \bigl|\bigl({\rm d}\bar L_v-{\rm d}\bar L_w\bigr)(\dot w,\xi,\dot w)\bigr|_{L^2}^2
\le c^2 \bigl(1+\abs{u}_{W^{1,2}} \bigr)^2
 \abs{v-w}_{W^{1,2}}^2
 \abs{\xi}_{W^{1,2}}^2\abs{w}_{W^{2,2}}^2,
\\
& \bigl|{\rm d}\bar L_v(\dot w,\xi,\dot v-\dot w)\bigr|_{L^2}^2
\le c^2 \bigl(1+\abs{u}_{W^{1,2}} \bigr)^2 \abs{\xi}_{W^{1,2}}^2
 \abs{v-w}_{W^{2,2}}^2,
\\
& \bigl|{\rm d}\bar L_v(\dot v-\dot w,\xi,\dot v)\bigr|_{L^2}^2
\le c^2 \bigl(1+\abs{u}_{W^{1,2}} \bigr)^2
 \abs{v-w}_{W^{2,2}}^2\abs{\xi}_{W^{1,2}}^2,
\\
& \bigl|\bigl(\bar L_v-\bar L_w\bigr) (\xi,\dot w)\bigr|_{L^2}^2
\le c^2 \abs{v-w}_{W^{1,2}}^2 \abs{\xi}_{W^{1,2}}^2 \abs{w}_{W^{2,2}}^2,
\\
& \bigl|\bar L_v(\xi,\dot v-\dot w)\bigr|_{L^2}^2
\le c^2 \abs{\xi}_{W^{1,2}}^2\abs{v-w}_{W^{2,2}}^2.
\end{split}
\end{equation*}
Continuing the above estimate
and returning to the notation $H_k=W^{k,2}$, we get
\begin{equation*}
\begin{split}
 \abs{(C^v-C^w)\xi}_{H_1}^2
 &\le \underbrace{7 c^2 \bigl(7+4\abs{u}_{H_1}^2\bigr)}_{=:\kappa^2}
 \bigl(\abs{v-w}_{H_2}^2+\abs{w}_{H_2}^2\abs{v-w}_{H_1}^2\bigr)
 \abs{\xi}_{H_1}^2 .
\end{split}
\end{equation*}
This implies that the operator norm is bounded by
\[
 \norm{C^v-C^w}_{\Ll(H_1)}
 \le\kappa (\abs{v-w}_{H_2}+\abs{w}_{H_2}\abs{v-w}_{H_1} ).
\]
Interchanging the roles of $v$ and $w$, we get
\[
 \norm{C^v-C^w}_{\Ll(H_1)}
 \le\kappa (\abs{v-w}_{H_2}+\abs{v}_{H_2}\abs{v-w}_{H_1} ).
\]
The above two estimates imply the scale Lipschitz estimate
\begin{equation}\label{eq:C-paraDarboux}
 \norm{C^v-C^w}_{\Ll(H_1)}
 \le\kappa (\abs{v-w}_{H_2}
 +\min\{\abs{v}_{H_2},\abs{w}_{H_2}\}\abs{v-w}_{H_1}
 ),
\end{equation}
which is precisely~(\ref{eq:(C)}) in axiom~(\texttt{C}).
This proves Step~(\texttt{C}).
\end{proof}

\noindent
\textbf{Step~(\texttt{F}).} $\forall u\in U_1\colon$
$F^u_2:=B_u^{-1}\p_t |_{H_2}\colon H_2\to H_1$ is Fredholm of index zero.

\begin{proof}
The proof was given in~(\ref{eq:hjghjghj77799})
as a consequence of Theorem~\ref{thm:Fredholm-F_C}.
\end{proof}

The proof of Theorem~\ref{thm:paradarboux-Hess} is complete.
\end{proof}

\subsubsection*{Tool}

\begin{Lemma}\label{le:sum-quad}
For $k\in\N$ real numbers $a_1,\dots,a_k>0$, there is the inequality
\[
 \Biggl(\sum_{j=1}^k a_j\Biggr)^2
 \le k\sum_{j=1}^k a_j^2 .
\]
\end{Lemma}

\begin{proof}
Observe that
\[
 \Biggl(\sum_{j=1}^k a_j\Biggr)^2
 =\sum_{j=1}^k a_j^2
 +\sum_{1\le j<i\le k}\underbrace{2a_ja_i}_{\le a_j^2+a_i^2}
 \le(1+(k-1)) a_1^2+\dots+(1+(k-1)) a_k^2 .
\tag*{\qed}
\]
\renewcommand{\qed}{}
\end{proof}

\subsection{Perturbed case}

\begin{Lemma}\label{le:para-weak-pert}
The perturbed para-Darboux Hessians~\eqref{eq:Hess-op} defined by
\[
 A^u\colon\ H_1\to H_0 ,\qquad
 \eta\mapsto B_u^{-1}\dot \eta+\bar L_u(\eta,\dot u)-a^u\eta
\]
one for each $u\in U_1$,
determine a weak Hessian field $A$ on $U_1$.
\end{Lemma}

\begin{Theorem}\label{thm:paradarboux-Hess-pert}
The perturbed para-Darboux weak Hessian field $A$ is almost extendable.
Moreover, the pair $(F,C)$ defined for $u\in U_1$ by
\[
 F^u\colon\ \eta\mapsto B_u^{-1}\dot \eta - a^u\eta ,\qquad
 C^u\colon\ \eta\mapsto \bar L_u(\eta,\dot u)
\]
is a decomposition.
\end{Theorem}

To prove the lemma and the theorem, we need the following proposition.

\begin{Proposition}\label{prop:a^u-cont}
The map $u\mapsto a^u$ is element of $C^0(U_1,\Ll(H_1))$.
\end{Proposition}

To prove the proposition, we need the following lemma.

\begin{Lemma}\label{le:tech-cont}
Let $f\in C^1(\R,\R)$. Then the map
\[
 F(f)\colon\ W^{1,2}\bigl(\SS^1,\R\bigr)\to W^{1,2}\bigl(\SS^1,\R\bigr)
 ,\qquad
 u\mapsto f\circ u
\]
is well defined and continuous.
\end{Lemma}

\begin{proof}
The composition $f\circ u$ is $L^2$, because
since $u$ is $C^0$ the image $u\bigl(\SS^1\bigr)$ is compact,
hence the continuous map $f$ is $L^2$-integrable along the image.
The derivative $f^\prime|_u\cdot \dot u$ is a product of a~$C^0$-map
$f^\prime|_u$ and an $L^2$-map $\dot u$.
Such product is $L^2$. This shows well defined.

We prove continuity. Pick $\eps>0$.
Since $u\colon\SS^1\to\R$ is continuous, there exists a compact
interval $[a,b]$ containing the image of $u$.
Since $f$ and $f^\prime$ are continuous, there exists
\[
 0<\delta<\min\biggl\{1,\frac{\eps}{2\abs{f^\prime(u)}_{L^\infty}}\biggr\}
\]
such that for all $x,y\in [a-1,b+1]$ with $\abs{x-y}<\delta$ it holds
\[
 \abs{f(x)-f(y)}\le \frac{\eps}{\kappa}
 ,\qquad
 \bigl|f^\prime(x)-f^\prime(y)\bigr|\le \frac{\eps}{\kappa}
 ,\qquad
 \kappa:=\sqrt{10+8\abs{u}_{W^{1,2}}^2} .
\]
Now assume that $\abs{u-v}_{W^{1,2}}\le\delta$.
Since $\abs{u-v}_{C^0}\le \abs{u-v}_{W^{1,2}}\le\delta$
we have $\abs{u_t-v_t}\le\delta\le 1$ whenever $t\in\SS^1$.
Hence
\[
 \abs{f(u_t)-f(v_t)}\le \frac{\eps}{\kappa}
 ,\qquad
 \bigl|f^\prime(u_t)-f^\prime(v_t)\bigr|\le \frac{\eps}{\kappa} ,
\]
for every $t\in\SS^1$.
We estimate
\begin{equation*}
\begin{split}
 \abs{f\circ u-f\circ v}_{W^{1,2}}^2
 &=\abs{f\circ u-f\circ v}_{L^2}^2
 +\bigl|f^\prime(u)\cdot \dot u
{\color{gray}\,
-f^\prime(u)\cdot \dot v
+f^\prime(u)\cdot \dot v
}
 -f^\prime(v)\cdot \dot v\bigr|_{L^2}^2
\\
 &\le \abs{f\circ u-f\circ v}_{L^2}^2
 +2\bigl|f^\prime(u)\cdot \dot u-f^\prime(u)\cdot \dot v\bigr|_{L^2}^2
 +2\bigl|f^\prime(u)\cdot \dot v-f^\prime(v)\cdot \dot v\bigr|_{L^2}^2
\\
 &\le \frac{\eps^2}{\kappa^2}
 +2\bigl|f^\prime(u)\bigr|_{L^\infty}^2\abs{\dot u-\dot v}_{L^2}^2
 +2\frac{\eps^2}{\kappa^2}\abs{\dot v}_{L^2}^2
\\
 &\le \frac{\eps^2}{\kappa^2}
 +2\bigl|f^\prime(u)\bigr|_{L^\infty}^2
 \underbrace{\abs{u-v}_{W^{1,2}}^2}
 _{\le\delta^2\le\frac{\eps^2}{4\abs{f^\prime(u)}_{L^\infty}^2}}
 +4\frac{\eps^2}{\kappa^2}
 \bigl(
 \abs{u}_{W^{1,2}}^2+\underbrace{\abs{u-v}_{W^{1,2}}^2}_{\le\delta\le 1}
 \bigr)
\\
 &\le \frac{\eps^2}{\kappa^2}\bigl(5+4\abs{u}_{W^{1,2}}^2\bigr)
 +\frac{\eps^2}{2}=\eps^2 .
\end{split}
\end{equation*}
This proves continuity and concludes the proof of
Lemma~\ref{le:tech-cont}.
\end{proof}

\begin{proof}[Proof of Proposition~\ref{prop:a^u-cont}]
Since the function $h\colon\mathfrak{U}\to\R$ is of class $C^3$,
given $i,j=1,\dots,2n$, the matrix coefficients
$\p_i\p_j h\colon \mathfrak{U}\to\R$ of the Hessian are $C^1$ functions.
Hence, by Lemma~\ref{le:tech-cont}
and the fact that multiplication $W^{1,2}\times W^{1,2}\to W^{1,2}$ is
continuous, the proposition follows.
\end{proof}

\begin{proof}[Proof of Lemma~\ref{le:para-weak-pert}]
By Proposition~\ref{prop:a^u-cont},
the map $u\mapsto a^u$ lies in $C^0(U_1,\Ll(H_1))$,
hence in particular in
$C^0(U_1,\Ll(H_1,H_0)) \cap C^0(U_2,\Ll(H_2,H_1))$.
Since $A^u=A^u_0-a^u$ and $A^u_0$ is a weak Hessian field according to
Lemma~\ref{le:para-weak}, it follows that
$A$ is also element of $C^0(U_1,\Ll(H_1,H_0)) \cap C^0(U_2,\Ll(H_2,H_1))$.

The (\texttt{Symmetry}) axiom~(\ref{eq:u-0-symmetric}) holds true by
Lemmas~\ref{le:A^u_0-symm} and~\ref{le:a^u-symm}.
We show the (\texttt{Fredholm}) axiom.
Since the inclusions $H_1\INTO H_0$ and $H_2\INTO H_1$
are compact, it follows that for every $u\in U_1$
the operator $a^u|_{H_2}\colon H_2\to H_1$ is compact
and so is $a^u\colon H_1\to H_0$.
By Lemma~\ref{le:para-weak}, the operators $A^u_0$ and $A^u_0|_{H_2}$
are in particular Fredholm of index zero.
Since Fredholm property and index are stable under compact perturbation,
it follows that the same is true for $A^u$ and for $A^u_2$.
This proves Lemma~\ref{le:para-weak-pert}.
\end{proof}

\begin{proof}[Proof of Theorem~\ref{thm:paradarboux-Hess-pert}]
Since $\Ll(H_1)\subset \Ll(H_2,H_1)\cap\Ll(H_1,H_0)$,
it follows from Proposition~\ref{prop:a^u-cont}
that the map $u\mapsto A^u$ lies in
$\in C^0(U_1, \Ll(H_2,H_1)\cap\Ll(H_1,H_0))$.
Since the inclusions $H_1\INTO H_0$ and $H_2\INTO H_1$
are compact, it follows that for every $u\in U_1$
the operator $a^u|_{H_2}\colon H_2\to H_1$ is compact
and so is $a^u\colon H_1\to H_0$.

By the unperturbed case, Theorem~\ref{thm:paradarboux-Hess},
we know that \smash{$B_u^{-1}\p_t$} as a map $H_1\to H_0$ and as a~map $H_2\to
H_0$ is Fredholm of index zero;
this is a consequence of Theorem~\ref{thm:Fredholm-F_C}.
Since Fredholm property and index are stable under compact
perturbation, the same is true for $F^u$ and its restriction $F^u_2$.
In particular, $F$ satisfies axiom~(\texttt{F}) in Definition~\ref{def:alm.ext}.

That $C$ satisfies axiom~(\texttt{C}) was already shown in the proof
of Theorem~\ref{thm:paradarboux-Hess}.
This shows that~$A$ is almost extendable and $(F,C)$ is a
decomposition of $A$.
\end{proof}

\appendix

\section{Theorem of Rabier}

While the result of Rabier applies in greater generality,
we discuss it for Hilbert space pairs.
Let~$(H,W)$ be a \textit{Hilbert space pair},
that is $H$ and $W$ are Hilbert spaces
such that $W\subset H$, as sets, and
the inclusion map $\iota\colon W\to H$ is
a compact linear operator.
Let ${\rm i}\in\C$ be the imaginary unit element.
Consider a family $(A(s))_{s\in\R}$
of bounded linear operators $A(s)\colon W\to H$.
Assume that the following five conditions hold:
\begin{itemize}\setlength\itemsep{0ex}
\item[\rm (H1)]
 $(H,W)$ is a Hilbert space pair;
\item[\rm (H2)]
 the operator family $\R\to\Ll(W,H)$, $s\mapsto A(s)$, depends
 continuously on~$s$;
\item[\rm (H3)]
 there are invertible asymptotic limits $A_\mp\in \Ll(W,H)$
 in the sense that
 \[
 \lim_{s\to\mp\infty} \norm{A(s)-A_\mp}_{\Ll(W,H)}=0 ;
 \]
\item[\rm (H4)]
 $\forall s\in\R\cup\{\mp\infty\}$
 $\exists C_0(s),r_0(s)>0$ $\forall \alpha\in\R$:
 Whenever $\abs{\alpha}\ge r_0(s)$ the operator
 $ A(s)-{\rm i}\alpha \colon W\to H$ is an isomorphism and there is the
 estimate
 \[
 \bigl\|(A(s)-{\rm i}\alpha)^{-1}\bigr\|_{\Ll(H)}
 \le \frac{C_0(s)}{\abs{\alpha}} ;
 \]
\item[\rm (H5)]
 ${\rm i}\R\cap\spec A_\mp=\varnothing$.
\end{itemize}

\begin{Theorem}[Rabier~\cite{Rabier:2004a}]\label{thm:Rabier}
Under the assumptions~{\rm (H1)--(H5)} the operator
\[
 D_A:=\p_s+A
 \colon\ W^{1,p}_{H}\cap L^p_{W}\to L^p_H
\]
is Fredholm for every $p\in(1,\infty)$.
\end{Theorem}

\subsection{The unparametrized Rabier condition (H4)}

In this appendix, we do not deal with operator families
parametrized by $s\in\R$, but with
individual operators $A\in\Ll(W,H)$ on a Hilbert space pair $(H,W)$.
An exemption is Corollary~\ref{cor:symm-Fredholm}.

\begin{Definition}
A bounded linear operator $A\colon W\to H$ \textit{satisfies the
Rabier condition}~(H4), more precisely~\cite[equation~(1.4)]{Rabier:2004a}
in the case of Hilbert spaces, if there are constants $C_0,r_0>0$
such that $\forall \alpha\in\R$ with $\abs{\alpha}\ge r_0$ it holds that
\[
 {\rm i}\alpha\notin
 \spec A
 :=\{\lambda\in \C\mid\text{$A-\lambda\colon W\to H$ is not
 bijective}\}
\]
and
\begin{equation}\label{eq:Rabier}
 \bigl\|\alpha (A-{\rm i}\alpha)^{-1}\bigr|\_{\Ll(H)}
 \le C_0 .
\end{equation}
Here and throughout, we write ${\rm i}\alpha$ to abbreviate ${\rm i}\alpha\iota$.
\end{Definition}

\subsubsection*{Symmetric Fredholm operators of index zero}

\begin{Lemma}\label{le:Rabier-symmetric}
Suppose $A\in\Ll(W,H)$ is $(a)$ $H$-symmetric
and $(b)$ Fredholm of index zero.
Then for every real number $\alpha\not= 0$, there is the estimate
\[
 \bigl\|\alpha (A-{\rm i}\alpha)^{-1}\bigr\|_{\Ll(H)}
 \le 1 .
\]
In particular, $A$ satisfies the Rabier condition~{\rm (H4)}
for $C_0=1$ and every $r_0>0$.
\end{Lemma}

\begin{proof}
By (a) and (b), the spectrum of $A$ is real
and there exists an ONB $\{e_n\}_{n\in\N}$ of $H$
and $\{a_n\}_{n\in\N}\subset \R$ such that $Ae_n=a_ne_n$
for every $n\in\N$; see, e.g.,~\cite[Appendix~D]{Frauenfelder:2024d}.
Let $\alpha\in\R\setminus\{0\}$, then
\[
 \alpha(A-{\rm i}\alpha)^{-1} e_n
 =\frac{\alpha}{a_n-{\rm i}\alpha} e_n .
\]
Let $\xi\in H$ and write it in the form
$\xi=\sum_n\xi_ne_n$ for unique real numbers $\xi_n$.
Suppose that $\xi$ is of unit norm $\abs{\xi}_{H}^2=\sum_n\xi_n^2=1$.
Then
\[
 \alpha(A-{\rm i}\alpha)^{-1}\xi
 =\sum_n \frac{\alpha\xi_n}{a_n-{\rm i}\alpha} e_n
\]
and we estimate
\begin{equation*}
\begin{split}
 \bigl|\alpha(A-{\rm i}\alpha)^{-1}\xi\bigr|_{H}^2
 =\sum_n \frac{\alpha^2\xi_n^2}{\abs{a_n-{\rm i}\alpha}_{\C}^2}
 =\sum_n \frac{\alpha^2}{a_n^2+\alpha^2} \xi_n^2
 \le \sum_n \xi_n^2
 =1.
\end{split}
\end{equation*}
This proves Lemma~\ref{le:Rabier-symmetric}.
\end{proof}

\begin{Corollary}\label{cor:symm-Fredholm}
Let $(A(s) \colon W\to H)_{s\in\R}$ be a family of bounded linear operators
such that each family member $A(s)$ is
$H$-symmetric and Fredholm of index zero
and {\rm (H1)--(H3)} are satisfied.
Then the operator
\[
 D_A:=\p_s+A
 \colon\ W^{1,p}_{H}\cap L^p_{W}\to L^p_H
\]
is Fredholm for every $p\in(1,\infty)$.
\end{Corollary}

\begin{proof}
By Rabier's Theorem~\ref{thm:Rabier},
it remains to verify (H4) and (H5).
(H4) holds by Lemma~\ref{le:Rabier-symmetric}.
We prove~(H5).
Since $A(s)$ is $H$-symmetric for every $s\in\R$, the
same holds for the asymptotics~$A_\mp$.
Therefore, the spectrum of $A_\mp$ is real.
Since $A_\mp$ are invertible by (H3)
zero does not belong to the spectrum of $A_\mp$
and therefore ${\rm i}\R\cap\spec A_\mp=\varnothing$.
\end{proof}

\subsubsection*{Reformulation}

\begin{Lemma}\label{le:Rabier-v2}
Let $A\colon W\to H$ be a bounded linear operator.
Then $A$ satisfies the Rabier condition~{\rm (H4)} if and only if it
satisfies the following two conditions:
\begin{itemize}\itemsep=0ex
\item[$(i)$]
 $A$ is Fredholm of index zero.
\item[$(ii)$]
 There exists constants $C_0,r_0>0$ such that
 whenever $\abs{\alpha}\ge r_0$ it holds
 \begin{equation}\label{eq:Rabier-v2}
 \abs{(A-{\rm i}\alpha)\xi}_{H}
 \ge C_0\alpha \abs{\xi}_{H}
 \end{equation}
 for every $\xi\in H$.
\end{itemize}
\end{Lemma}

\begin{proof}
`$\Leftarrow$'
Since $A\colon W\to H$ is Fredholm of index zero and inclusion
$W\INTO H$ is compact,
it follows that $A-{\rm i}\alpha\colon W\to H$ is also Fredholm of index
zero for every $\alpha\in\R$.
If $r_0$ is as in (ii) and $\alpha\ge r_0$, then $A-{\rm i}\alpha$ is
injective by the estimate in (ii), hence surjective (Fredholm index zero).
Therefore, by the open mapping theorem, the operator $A-{\rm i}\alpha$ is
an isomorphism. In particular, the inverse \smash{$(A-{\rm i}\alpha)^{-1}$} is bounded.
Now~(\ref{eq:Rabier}) follows from~(\ref{eq:Rabier-v2}).

`$\Rightarrow$'
Suppose (H4) holds true.
Hence, if $\alpha\ge r_0$, then $A-{\rm i}\alpha$ is bijective
and therefore it is Fredholm of index zero.
Since inclusion $W\INTO H$ is compact, it follows that $A$
is Fredholm of index zero as well.
This shows (i) and~(\ref{eq:Rabier-v2}) follows from~(\ref{eq:Rabier}).
This proves Lemma~\ref{le:Rabier-v2}.
\end{proof}

\begin{Remark}[spectrum consists of eigenvalues]
$\lambda\in\spec A$ is called
\textit{eigenvalue} of $A\in\Ll(W,H)$
if $A-\lambda\iota\colon W\to H$ is not injective.
Observe that
\begin{equation*}
 \text{$A$ satisfies (H4)}\quad\Rightarrow\quad
 \spec A=\{\text{eigenvalues of $A$}\}.
\end{equation*}
To see this, suppose $A$ satisfies~(H4) and pick $\lambda\in \spec A$.
The latter means that $A-\lambda\iota$ is not injective or not surjective.
Now $A$ is Fredholm of index zero by~(i) and so is $A-\lambda\iota$,
as shown in the beginning of the previous proof.
In particular, not injective and not surjective are equivalent for
the operator $A-\lambda\iota$.
\end{Remark}

\subsubsection*{Stability under compact perturbation}

The following proposition is a version of~\cite[Theorem~3.5]{Rabier:2003a}.

\begin{Proposition}[compact perturbation]\label{prop:Rabier+compact}
Let $A\colon W\to H$ be a bounded linear operator
satisfying the Rabier condition~{\rm (H4)}
and $K\colon W\to H$ a compact linear operator.
Then $A+K\colon W\to H$
satisfies the Rabier condition~{\rm (H4)} as well.
\end{Proposition}

In the proof, we utilize the following notions.
Let $A\colon W\to H$ be a bounded linear operator.
A bounded linear operator $T\colon W\to H$ is called
\textit{\boldmath$A$-bounded} if
there exist constants $a\ge 0$ and $b(a)\ge 0$ such that
\begin{equation}\label{eq:A-bounded}
 \abs{T\xi}_H\le a\abs{A\xi}_H+b(a)\abs{\xi}_H
\end{equation}
for every $\xi\in W$.
The \textit{\boldmath$A$-bound} of $T$
is the infimum of all possible values $a\ge 0$, in symbols
\[
 \inf\{a\ge 0\mid\text{$\exists b(a)$ such that~(\ref{eq:A-bounded})
 holds $\forall \xi\in W$}\}.
\]
A theorem of Hess~\cite[Satz~1]{Hess:1969a}
says that, given $A\in\Ll(W,H)$ with non-empty resolvent set,
compact operators $K\colon W\to H$ have vanishing $A$-bound.

\begin{proof}
Let $A$ satisfy~(H4). Then $A$ satisfies (i) and (ii) in
Lemma~\ref{le:Rabier-v2} with constants $r_0$, $C_0$.
It~suffices to check conditions (i) and (ii)
for $A+K$ for some constants~$r_1$,~$C_1$.
Since $K\colon W\to H$ is compact,
it holds that $A+K$ is Fredholm and
$\INDEX (A+K)=\INDEX A=0$. This proves (i).
To prove (ii) we define
\[
 \eps:=\min\biggl\{\frac12,\frac{C_0}{4}\biggr\}
 ,\qquad
 r_1:=\max\biggl\{r_0,\frac12+4\frac{b(\eps)}{C_0}\biggr\}
 ,\qquad
 C_1:=\frac{C_0}{8}
\]
and for $\abs{\alpha}\ge r_1\ge r_0$ we estimate
\begin{align*}
\begin{split}
 \abs{(A+K-{\rm i}\alpha)\xi}_{H}
 &\ge\abs{(A-{\rm i}\alpha)\xi}_{H}-\abs{K\xi}_{H}
\\
 &\ge\abs{(A-{\rm i}\alpha)\xi}_{H}
 -\eps \abs{(A{\color{gray}\,-{\rm i}\alpha+{\rm i}\alpha})\xi}_{H}
 -b(\eps) \abs{\xi}_{H}
\\
 &\ge (1-\eps) \abs{(A-{\rm i}\alpha)\xi}_{H}
 - (\eps\abs{\alpha}+b(\eps) )\abs{\xi}_H
\\
 &\ge ((1-\eps) C_0\abs{\alpha}-\eps\abs{\alpha}-b(\eps) )
 \abs{\xi}_H
\\
 &\ge\biggl(\frac{C_0}{4}\abs{\alpha}-b(\eps)\biggr)
 \abs{\xi}_H
\\
 &\ge\frac{C_0}{8}\abs{\xi}_H .
 \end{split}
\end{align*}
In inequality~2, we used~(\ref{eq:A-bounded}) and then added zero
$-{\rm i}\alpha+{\rm i}\alpha$.
Inequality~3 is by the triangle inequality.
Inequality~4 is by the hypothesis~(\ref{eq:Rabier-v2}) on $A$
which applies since $\abs{\alpha}\ge r_0$.
Inequality~5 is by choice of $\eps$.
Inequality~6 uses $\abs{\alpha}\ge r_1$
and the choice of $r_1$.
This proves Proposition~\ref{prop:Rabier+compact}.
\end{proof}

\subsubsection*{Level operator}

\begin{Lemma}\label{le:Rabier+level-2}
Let $(H_0,H_1,H_2)$ be a Hilbert space triple.
Suppose $F\colon H_1\to H_0$ is a bounded linear operator
which satisfies the Rabier estimate~{\rm (H4)}.
If $F$ restricts to a map ${F_2\!=\!F|_{H_2}\!\colon
H_2\!\to\! H_1}$ and $F_2$ is Fredholm of index zero,
then $F_2$ satisfies the Rabier estimate~{\rm (H4)} as well.
\end{Lemma}

\begin{proof}
Let $C_0,r_0>0$ be the constants in the Rabier condition~(H4).
Fix $\alpha_0\ge r_0$.
Then $F-{\rm i}\alpha_0\colon H_1\to H_0$ is bijective, thus an isomorphism.
Therefore, maybe after replacing the norm of $H_1$ by an equivalent
norm, we can assume without loss of generality
that for every $\xi\in H_1$ the $H_1$-norm is given by
$\abs{\xi}_1 =\abs{(F-{\rm i}\alpha_0)\xi}_0$.
Hence for $\alpha\ge r_0$ and $\xi\in H_1$, we write and estimate
\begin{equation*}
\begin{split}
 \abs{(F-{\rm i}\alpha)\xi}_1
 &=\abs{(F-{\rm i}\alpha_0) (F-{\rm i}\alpha)\xi}_0\\
 &=\abs{(F-{\rm i}\alpha) (F-{\rm i}\alpha_0)\xi}_0\\
 &\ge C_0\alpha \abs{(F-{\rm i}\alpha_0)\xi}_0\\
 &= C_0\alpha \abs{\xi}_1 .
\end{split}
\end{equation*}
The inequality is by~(\ref{eq:Rabier-v2}) for $F$.
This shows that $F_2$ satisfies~(\ref{eq:Rabier-v2}).
But $F_2$ is Fredholm of index zero by hypothesis.
Hence $F_2$ satisfies (i) and (ii) in
Lemma~\ref{le:Rabier-v2} and therefore~(H4).
This proves Lemma~\ref{le:Rabier+level-2}.
\end{proof}

\section{Compact non-continuous perturbations}

\begin{Lemma}\label{le:bhgb777}
For a Hilbert space $H$,
an interval $I=[-T,T]$, and $c\in L^2(I,H)$
multiplication $m_c\colon W^{1,2}(I,\R)\to L^2(I,H)$,
$\xi\mapsto c\xi$ is a compact operator.
\end{Lemma}

\begin{proof}
The multiplication operator is a composition
\[
 m_c=m_c^0\circ \iota\colon\
 W^{1,2}(I,\R)
 \stackrel{\iota}{\longrightarrow} C^0(I,\R)
 \stackrel{m_c^0}{\longrightarrow} L^2(I,H)
\]
of the compact Sobolev inclusion $\iota\colon W^{1,2}(I,\R)\INTO C^0(I,\R)$
followed by the multiplication $m_c^0$
which is bounded as we show next.
Indeed, multiplication
\[
 L^2(I,H)\times C^0(I,\R)\to L^2(I,H)
 ,\qquad
 (c,\xi)\mapsto c\xi
\]
is a continuous map since
\[
 \bigl\|m_c^0\xi\bigr\|_{L^2(I,H)}
 =\norm{c\xi}_{L^2(I,H)}
 \le \norm{c}_{L^2(I,H)}\norm{\xi}_{C^0(I,\R)} .
\]
Hence the map $m_c$
is the composition of a compact and a continuous
linear map, and therefore it is compact.
This proves Lemma~\ref{le:bhgb777}.
\end{proof}

\begin{Proposition}\label{prop:bhgb777}
Assume $H$ is a Hilbert space and $c\in L^2(\R,H)$.
Then multiplication $m_c\colon W^{1,2}(\R,\R)\to L^2(\R,H)$,
$\xi\mapsto c\xi$ is a compact operator.
\end{Proposition}

\begin{proof}
Given $n\in\N$, let $\chi_n:=\chi_{[-n,n]}$ be the characteristic
function on the interval $[-n,n]$, that is $\chi_n$ is $1$ on $[-n,n]$
and $0$ else. Then the multiplication operator $m^n_c:=m_{\chi_n c}$
is compact by Lemma~\ref{le:bhgb777}.

To show that $m_c$ is compact, we show that the sequence of compact
operators $m^n_c$ converges in the operator norm topology to $m_c$.
Indeed,
\begin{equation*}
\begin{split}
 \norm{(m_c-m^n_c)\xi}_{L^2(\R,H)}
 &=\norm{(c-\chi_n c)\xi}_{L^2(\R,H)}\\
 &=\bigl\|\chi_{\R\setminus [-n,n]} c\xi\bigr\|_{L^2(\R,H)}\\
 &\le \bigl\|\chi_{\R\setminus [-n,n]} c\bigr\|_{L^2(\R,H)}
 \norm{\xi}_{C^0(\R,\R)}\\
 &\le \underbrace{\norm{c}_{L^2(\R\setminus[-n,n],H)}}
 _{\text{$\to 0$, as $n\to\infty$}}
 \norm{\xi}_{W^{1,2}(\R,\R)} .
\end{split}
\end{equation*}
This proves Proposition~\ref{prop:bhgb777}.
\end{proof}

Let $(H_1,H_2)$ be a \textit{Hilbert space pair},
see, e.g.,~\cite{Frauenfelder:2024c},
that is $H_1$ and $H_2$ are both infinite-dimensional
Hilbert spaces such that $H_2\subset H_1$ and inclusion
$\iota\colon H_2\INTO H_1$ is compact and dense.
We define \smash{$W^{1,2}_{H_1}$} and \smash{$L^2_{H_i}$} by~(\ref{eq:abbreviate}).

\begin{Theorem}\label{thm:perturbation}
Let $C\in L^2(\R,\Ll(H_1))$.
Then pointwise multiplication
\[
 M_C\colon\
 W^{1,2}_{H_1}\cap L^2_{H_2}\to L^2_{H_1}
 ,\qquad \xi\mapsto [s\mapsto C(s)\xi(s) ]
\]
is a compact operator.
\end{Theorem}

\begin{Remark}
In the special case where $C$ is continuous,
Theorem~\ref{thm:perturbation}
was proved by Robbin and Salamon~\cite[Lemma~3.18]{robbin:1995a}.
As in our proof, they
use a sequence of compact operators which converges
to the operator $M_C$ in the operator topology.
In contrast to our sequence, they first multiply with~$C$ and
then project, while we first project and then multiply by~$C$.
\end{Remark}

\begin{Definition}[intersection of Hilbert spaces]\label{def:intersection}
The intersection $H\cap W$ of two Hilbert spaces~$H$ and~$W$ is itself
a~Hilbert space with inner product and norm
\[
 \INNER{\cdot}{\cdot}_{H\cap W}
 :=\INNER{\cdot}{\cdot}_H+\INNER{\cdot}{\cdot}_W
 ,
 \qquad
 \norm{\cdot}_{H\cap W}
 :=\sqrt{\norm{\cdot}_H^2+\norm{\cdot}_W^2} .
\]
\end{Definition}

\begin{proof}[Proof of Theorem~\ref{thm:perturbation}]
The Hilbert space pair $(H_1,H_2)$ is isometric to the pair
$\bigl(\ell^2,\ell^2_h\bigr)$ where $h\colon\N\to(0,\infty)$ is a monotone
unbounded function; see~\cite[Theorem~A.4]{Frauenfelder:2024c}.
In the following, we identify the Hilbert space pairs
\[
 (H_1,H_2)\simeq \bigl(\ell^2,\ell^2_h\bigr).
\]
We denote by $\Ee=\{e_k\}_{k\in\N}$ the canonical basis of $\ell^2$
and by
\[
 \pi_n\colon\ H_1\to\R^n
 ,\qquad
 \xi=\sum_{k=1}^\infty \xi_k e_k
 \mapsto \left(\xi_1,\dots,\xi_n\right)
\]
the orthogonal projection of $H_1$ to the $n$-dimensional subspace of
$H_1$ identified with $\R^n$ under the isometry mentioned above.
As the basis $\Ee$ is still orthogonal in \smash{$\ell^2_h$},
the restriction ${\pi_n|_{H_2}\colon H_2\to\R^n}$ is again the
orthogonal projection.

For $n\in\N$, we define the operator
\begin{equation*}
\begin{split}
 M^n_C\colon\ W^{1,2}_{H_1}\cap L^2_{H_2}&\to L^2_{H_1}
 ,\qquad
 \xi\mapsto [s\mapsto
 (M^n_C \xi )(s):=C(s)\pi_n \xi(s) ]
\end{split}
\end{equation*}
where $\iota_n=\pi_n^*\colon\R^n\INTO H_1$ is inclusion.

\medskip

 \noindent
 \textbf{Step~1.}
For every $n\in\N$, the operator $M^n_C$ is compact.

Step~1 follows from Proposition~\ref{prop:bhgb777}
for $m_c=M^n_C$ and where \smash{$\tilde\xi(s)=\pi_n\xi(s)$} now takes values
in~$\R^n$ instead of $\R$.

\medskip

 \noindent
\textbf{Step~2.}
As $n\to\infty$, the operators $M^n_C$
converge to $M_C$ in the operator norm
\[
 \lim_{n\to\infty}
 \sup_{\norm{\xi}_{W^{1,2}_{H_1}\cap L^2_{H_2}=1}}
 \Norm{(M^n_C-M_C)\xi}_{L^2_{H_1}}
 =0 .
\]

\medskip \noindent
To prove Step~2, we write the difference in the form
\[
 M_C-M_C^n=\Mm_C\circ P_n,
\]
where
\[
 P_n\colon\
 W^{1,2}_{H_1}\cap L^2_{H_2}\to L^\infty_{H_1}
 ,\qquad \xi\mapsto [s\mapsto (\1-\pi_n )\xi(s) ]
\]
and
\[
 \Mm_C\colon\
 L^\infty_{H_1}\to L^2_{H_1}
 ,\qquad \eta\mapsto [s\mapsto C(s)\eta(s) ] .
\]
The operator $\Mm_C$ is bounded: Indeed, for
$\eta\in L^\infty_{H_1}$ of unit norm, we estimate
\begin{equation*}
\begin{split}
 \norm{\Mm_C \eta}_{L^2_{H_1}}^2
 &=\int_{-\infty}^\infty\abs{C(s)\eta(s)}_{H_1}^2\, {\rm d}s\\
 &\le \int_{-\infty}^\infty\norm{C(s)}_{\Ll(H_1)}^2
 \Bigl(\sup_{s\in\R}\abs{\eta(s)}_{H_1} \Bigr)^2\, {\rm d}s\\
 &=\norm{C}_{L^2(\R,\Ll(H_1))}^2 .
\end{split}
\end{equation*}
Since
\begin{equation*}
\begin{split}
 \norm{M_C-M^n_C}_{\Ll(W^{1,2}_{H_1}\cap L^2_{H_2},L^2_{H_1})}
 \le
 \norm{\Mm_C}_{\Ll(L^\infty_{H_1}, L^2_{H_1})}\cdot
 \norm{P_n}_{\Ll(W^{1,2}_{H_1}\cap L^2_{H_2},L^\infty_{H_1})},
\end{split}
\end{equation*}
it suffices to show that the norm of $P_n$ converges to zero, as
$n\to\infty$.

\begin{Claim*}
If \smash{$\eta\in W^{1,2}_{H_1}\cap L^2_{H_2}$} satisfies the conditions
\[
 \text{\rm (i) $\norm{\dot \eta}_{L^2_{H_1}}\le1$}
 ,\qquad
 \text{\rm (ii) $\norm{\eta}_{L^2_{H_2}}\le1$}
 ,\qquad
 \text{\rm (iii) $\forall s\in\R\colon\;\pi_n\eta(s)=0$,}
\]
then
\smash{$
 \norm{\eta}_{L^\infty_{H_1}}
 \le\bigl(\frac{3}{h_{n+1}}\bigr)^\frac{1}{4}
$}.
\end{Claim*}

\begin{proof}
Condition (iii) means that $\eta$ is for any $s\in\R$ of the form
\[
 \eta(s)
 \stackrel{\text{(iii)}}{=}
 (0,\dots,0,\eta_{n+1}(s),\eta_{n+2}(s),\dots )
\]
which, as in~(\ref{eq:spike-f}), provides the estimate
\begin{equation}\label{eq:bkjh7447ggg7}
 \abs{\eta(s)}_{\ell^2_h}^2
 \ge h_{n+1}\abs{\eta(s)}_{\ell^2}^2 .
\end{equation}
Now we prove the estimate
\begin{equation}\label{eq:claim-s}
 \abs{\eta(s)}_{\ell^2}
 \ge \abs{\eta(0)}_{\ell^2}-\sqrt{\abs{s}}
\end{equation}
for every $s\in\R$.
We show this estimate for $s\ge0$, in case $s\le0$ the argument is similar.
To this end pick $s\ge0$. Use Cauchy--Schwarz inequality to estimate
\[
 \int_0^s 1\cdot \abs{\dot \eta(t)}_{\ell^2}\, {\rm d}t
 \le\sqrt{\int_0^s 1\, {\rm d}t}\cdot
 \sqrt{\int_0^s \abs{\dot \eta(t)}_{\ell^2}\, {\rm d}t}
 \le \sqrt{s}\cdot\norm{\dot\eta}_{L^2_{H_1}}
 \stackrel{\text{(i)}}{\le}\sqrt{s} .
\]
This estimate and the fundamental theorem of calculus
yield the claim, namely
\begin{gather*}
 \abs{\eta(s)}_{\ell^2}
 =\Abs{\eta(0)+\int_0^s\dot\eta(t)\, {\rm d}t}_{\ell^2}
 \ge \abs{\eta(0)}_{\ell^2}
 -\int_0^s 1\cdot \abs{\dot \eta(t)}_{\ell^2}\, {\rm d}t
 \ge \abs{\eta(0)}_{\ell^2}-\sqrt{s} .
\end{gather*}
This proves the estimate~(\ref{eq:claim-s}).
Next we estimate
\begin{equation*}
\begin{split}
 1
 &\stackrel{\text{(ii)}}{\ge}\norm{\eta}_{L^2_{H_2}}^2\\
 &=\int_{\infty}^\infty\abs{\eta(s)}_{\ell^2_h}^2\, {\rm d}s\\
 &\stackrel{3}{\ge} h_{n+1}\int_{\infty}^\infty\abs{\eta(s)}_{\ell^2}^2\, {\rm d}s\\
 &\ge h_{n+1}\int_{-\abs{\eta(0)}_{\ell^2}^2}^{\abs{\eta(0)}_{\ell^2}^2}
 \abs{\eta(s)}_{\ell^2}^2\, {\rm d}s\\
 &\stackrel{5}{\ge} 2h_{n+1}\int_0^{\abs{\eta(0)}_{\ell^2}^2}
 \bigl(\abs{\eta(0)}_{\ell^2}-\sqrt{s} \bigr)^2 {\rm d}s\\
 &=\frac{1}{3} h_{n+1} \abs{\eta(0)}_{\ell^2}^4.
\end{split}
\end{equation*}
In Step~3, we used~(\ref{eq:bkjh7447ggg7}).
In Step~5, we used the claim~(\ref{eq:claim-s}).
To obtain the final identity, we calculated the integral.

We rewrite the obtained estimate in the form
\smash{$\abs{\eta(0)}_{\ell^2}\le (3/h_{n+1})^{\frac14}$}.
For $r\in\R$, we define $\eta_r(s):=\eta(s+r)$.
Then conditions (i)--(iii) hold as well for $\eta_r$.
Therefore,
\[
\abs{\eta(r)}_{\ell^2}=\abs{\eta_r(0)}_{\ell^2}
\le (3/h_{n+1})^{\frac14}
\] for every $r\in\R$.
Consequently,
\smash{$
 \norm{\eta}_{L^\infty_{H_1}}
 \le (3/h_{n+1} )^\frac{1}{4}
$}.
This proves the claim.
\end{proof}

\textit{Operator norm of \boldmath$P_n$.}
Pick \smash{$\xi\in W^{1,2}_{H_1}\cap L^2_{H_2}$} such that
\smash{$\norm{\xi}_{W^{1,2}_{H_1}\cap L^2_{H_2}}\le1$}.
Now we verify conditions (i)--(iii) for $\eta=P_n\xi$.
Condition~(i) is satisfied, indeed
\begin{gather*}
\begin{split}
 \bigl\|P_n\dot\xi\bigr\|_{L^2_{H_1}}
 &=\bigl\|(\1-\pi_n)\dot\xi\bigr\|_{L^2_{H_1}}
 \le\bigl\|\dot\xi\bigr\|_{L^2_{H_1}}
 \le\norm{\xi}_{W^{1,2}_{H_1}}
 \le\norm{\xi}_{W^{1,2}_{H_1} \cap L^2_{H_2}}
 \le 1 .
\end{split}
\end{gather*}
Inequality one uses that the projection is $H_1$-orthogonal.
In inequality three, we used Definition~\ref{def:intersection}
of the norm in an intersection space.
Condition~(ii) is satisfied by the same arguments
\begin{equation*}
\begin{split}
 \norm{P_n\xi}_{L^2_{H_2}}
 =\norm{(\1-\pi_n)\xi}_{L^2_{H_2}}
 \le\norm{\xi}_{L^2_{H_2}}
 \le 1 .
\end{split}
\end{equation*}
In inequality one, we used that the projection is also $H_2$-orthogonal.
Condition~(iii) is satisfied due to the projection property
$\pi_n^2=\pi_n$, namely
\[
 \pi_nP_n\xi(s)
 =\pi_n(\1-\pi_n)\xi(s)
 =0 .
\]
Thus, by the claim, we get
\smash{$
 \norm{P_n\xi}_{L^\infty_{H_1}}
 \le\bigl(\frac{3}{h_{n+1}}\bigr)^\frac{1}{4}
$}.
Therefore, the operator norm
\[
 \Norm{P_n}_{\Ll(W^{1,2}_{H_1}\cap L^2_{H_2},L^\infty_{H_1})}
 \le \biggl(\frac{3}{h_{n+1}}\biggr)^\frac{1}{4}
\]
converges to zero, as $n\to\infty$, since the growth function $h$ is
unbounded.
This concludes the proof of Step~2.

Steps~1 and~2 together imply Theorem~\ref{thm:perturbation}
by the standard fact that the subspace of compact operators
is closed in the space of bounded linear operators;
see, e.g.,~\cite[Theorem~A.1]{Frauenfelder:2021b}.
This proves Theorem~\ref{thm:perturbation}.
\end{proof}

\begin{Example}\label{ex:spike-f}
Recall that we identify
\smash{$(H_1,H_2)\simeq \bigl(\ell^2,\ell^2_h\bigr)$}.
Pick an element $\eta\in\ell^2$ which is of the form
$\eta=(0,\dots,0,\eta_{n+1},\eta_{n+2},\dots)$.
Set $a:=\abs{\eta}_{\ell^2}$ and pick $b>0$.
Define a map $\xi\colon\R\to\ell^2$~by\looseness=-1
\begin{equation*}
 \xi(s):=
 \begin{cases}
 \bigl(\frac{s}{b}+1\bigr)\eta,&\text{$s\in[-b,0]$,}
 \\
 \bigl(-\frac{s}{b}+1\bigr)\eta,&\text{$s\in[0,b]$,}
 \\
 0,&\text{else,}
 \end{cases}
\end{equation*}
as illustrated by Figure~\ref{fig:fig-spike-f}.

\begin{figure}[t]
 \centering
 \includegraphics
 {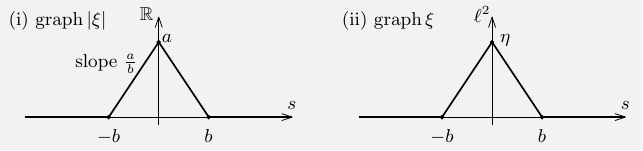}
 \caption{Graphs of $\abs{\xi}$ and $\xi$ in Example~\ref{ex:spike-f}.}
 \label{fig:fig-spike-f}
\end{figure}

Note that $\abs{\dot\xi(s)}_{\ell^2}=a/b$
whenever $s\in(-b,0)\cup(0,b)$.
Suppose that $\xi$ satisfies the conditions
(a) \smash{$\bigl\|\dot \xi\bigr\|_{L^2_{H_1}}=1$}
and
(b) $\norm{\xi}_{L^2_{H_2}}=1$.
By (a) and calculation,
\[
 1
 =\bigl\|\dot \xi\bigr\|_{L^2_{H_1}}^2
 =2\int_0^b \Bigl(\frac{a}{b}\Bigr)^2 {\rm d}s
 =\frac{2a^2}{b}
 \quad\Rightarrow\quad
 b=2a^2
\]
and
\[
 \norm{\xi}_{L^2_{H_1}}^2
 =2\int_0^b \Bigl(\frac{a}{b} s\Bigr)^2 {\rm d}s
 =2a^2\frac{b}{3}
 =\frac{4}{3} a^4 .
\]
Observe that
\begin{equation}\label{eq:spike-f}
 \abs{\eta}_{\ell^2_h}^2
 =\sum_{k=n+1}^\infty h_k\eta_k^2
 \ge h_{n+1} \sum_{k=n+1}^\infty \eta_k^2
 =h_{n+1}\abs{\eta}_{\ell^2}^2
\end{equation}
and use the unit norm condition (b) to obtain
\[
 1=\norm{\xi}_{L^2_{H_2}}^2
 \ge h_{n+1} \norm{\xi}_{L^2_{H_1}}^2
 =\frac{4}{3} h_{n+1} a^4,
\]
which implies that
$
 a\le\bigl(\frac{3}{4h_{n+1}}\bigr)^\frac{1}{4}
 \longrightarrow 0
$,
as $n\to\infty$.
\end{Example}

\section{Slice-wise index zero Fredholm operators}\label{sec:Fred-ind}

The space of symplectic bilinear forms on $\R^{2n}$ is denoted by
\[
 \Omega\bigl(\R^{2n}\bigr)
 =\bigl\{\text{$\omega\colon \R^{2n}\times \R^{2n}\to \R$
 bilinear, skew-symmetric, non-degenerate}\bigr\} .
\]
A loop of symplectic forms
\begin{equation}\label{eq:loop-omega}
 \omega\in W^{1,2}\bigl(\SS^1,\Omega\bigl(\R^{2n}\bigr)\bigr)
 ,\qquad
 t\mapsto \omega_t
\end{equation}
determines a loop of invertible matrices
$B_\omega\in W^{1,2}\bigl(\SS^1,\GL\bigl(2n,\R\bigr)\bigr)$ by
\[
 \omega_t (\cdot,B_\omega(t)\cdot )=\INNER{\cdot}{\cdot}
\]
pointwise for every $t\in\SS^1$ where $\INNER{\cdot}{\cdot}$ is the
Euclidean inner product on $\R^{2n}$.

\begin{Theorem}\label{thm:Fredholm-F_C}
Fix a loop of symplectic bilinear forms
\smash{$\omega\in W^{1,2}\bigl(\SS^1,\Omega\bigl(\R^{2n}\bigr)\bigr)$}.
Then the operator
\[
 F_{B_\omega^{-1}}=B_\omega^{-1}\p_t
 \colon\ H_{k}\to H_{k-1}
 ,\qquad
 H_k:=W^{k,2}\bigl(\SS^1,\R^{2n}\bigr) ,
\]
is Fredholm of index zero
whenever $k=1,2$.
\end{Theorem}

The proof of Theorem~\ref{thm:Fredholm-F_C}
will be given at the end of Appendix~\ref{sec:Fred-ind}.

\subsubsection*{Semi-Fredholm estimates}

\begin{Lemma}\label{le:sF-1}
Consider two loops of invertible matrices
\[
 C\in C^0\bigl(\SS^1,\GL(k,\R)\bigr)
 ,\qquad
 L\in W^{1,2}\bigl(\SS^1,\GL(k,\R)\bigr) .
\]
Then the corresponding linear operators
\begin{equation}\label{eq:sF-1}
\begin{split}
 \text{$(a)$}\
 F_C\colon\ W^{1,2}\bigl(\SS^1,\R^k\bigr)&\to L^2 \bigl(\SS^1,\R^k\bigr)
 ,\qquad
 \xi\mapsto \bigl[t\mapsto C(t)\dot\xi(t)\bigr]
\\
 \text{$(b)$}\
 F_L\colon\ W^{2,2}\bigl(\SS^1,\R^k\bigr)&\to W^{1,2}\bigl(\SS^1,\R^k\bigr)
 ,\qquad
 \xi\mapsto \bigl[t\mapsto L(t)\dot\xi(t)\bigr]
\end{split}
\end{equation}
are semi-Fredholm operators.\footnote{Semi-Fredholm means finite-dimensional kernel and closed range.}
\end{Lemma}

\begin{proof}
(a)~Since $C(t)$ is invertible and the circle is compact, there exists
$\delta>0$ such that $\abs{C(t)v}\ge \delta \abs{v}$
for all $t\in\SS^1$ and $v\in\R^k$.
Hence we estimate
\begin{equation*}
\begin{split}
 \abs{F_C\xi}_{L^2}^2
 =\int_0^1\bigl|C(t)\dot\xi(t)\bigr|^2\, {\rm d}t
 \ge\delta^2 \bigl|\dot\xi\bigr|_{L^2}^2
 =\delta^2 \abs{\xi}_{W^{1,2}}^2-\delta^2\abs{\xi}_{L^2}^2 .
\end{split}
\end{equation*}
Therefore,
\[
 \abs{\xi}_{W^{1,2}}^2
 \le\frac{1}{\delta^2}\abs{F_C\xi}_{L^2}^2+\abs{\iota \xi}_{L^2}^2
\]
for every $\xi\in W^{1,2}$.
Since the inclusion map $\iota\colon W^{1,2}\to L^2$ is compact,
this is a semi-Fredholm estimate, see, e.g.,~\cite[Lemma~A.1.1]{mcduff:2004a}.
It implies finite dimension of the kernel and closedness of the image
of the operator $F_C$.

(b)~Since $W^{1,2}$ embeds in $C^0$, the map $L$ is continuous.
Since $L(t)$ is invertible and the circle is compact, there exists
$\delta>0$ such that $\abs{L(t)v}\ge \delta \abs{v}$
for all $t\in\SS^1$ and $v\in\R^k$.
Hence we estimate
\begin{equation*}
\begin{split}
 \abs{F_L\xi}_{W^{1,2}}^2
 &=\bigl|L\dot\xi\bigr|_{L^2}^2+\bigl|\dot L\dot\xi+L\ddot\xi\bigr|_{L^2}^2
\\
 &\ge\delta^2 \bigl|\dot\xi\bigr|_{L^2}^2
 +\bigl|\dot L\dot\xi\bigr|_{L^2}^2
 \underline{+2\biggl\langle\!{\sqrt{2}\dot L\dot\xi},{\frac{1}{\sqrt{2}}L\ddot\xi}\biggr\rangle}
 +\bigl|L\ddot\xi\bigr|_{L^2}^2
\\
 &\ge\delta^2 \bigl|\dot\xi\bigr|_{L^2}^2
 +\bigl|\dot L\dot\xi\bigr|_{L^2}^2
 \underline{-2\bigl|\dot L\dot\xi\bigr|_{L^2}^2-\frac12\bigl|L\ddot\xi\bigr|_{L^2}^2}
 +\bigl|L\ddot\xi\bigr|_{L^2}^2
\\
 &=\delta^2 \bigl|\dot\xi\bigr|_{L^2}^2
 +\frac12 \bigl|L\ddot\xi\bigr|_{L^2}^2
 -\bigl|\dot L\dot\xi\bigr|_{L^2}^2
\\
 &\ge \delta^2 \bigl|\dot\xi\bigr|_{L^2}^2
 +\frac{\delta^2}{2}\bigl|\ddot\xi\bigr|_{L^2}^2
 -\bigl|\dot L\bigr|_{L^2}^2 \bigl|\dot\xi\bigr|_{C^0}^2
\\
 &\ge \frac{\delta^2}{2}\abs{\xi}_{W^{2,2}}^2
 -\frac{\delta^2}{2}\abs{\xi}_{L^2}^2
 -\abs{\dot L}_{L^2}^2 \abs{\xi}_{C^1}^2
\\
 &\ge \frac{\delta^2}{2}\abs{\xi}_{W^{2,2}}^2
 -\biggl(\frac{\delta^2}{2}+\abs{\dot L}_{L^2}^2 \biggr) \abs{\xi}_{C^1}^2.
\end{split}
\end{equation*}
Therefore,
\[
 \abs{\xi}_{W^{2,2}}^2
 \le\frac{2}{\delta^2}\abs{F_L\xi}_{L^2}^2
 +\biggl(1+\frac{2}{\delta^2}\abs{\dot L}_{L^2}^2\biggr)
 \abs{\xi}_{C^1}^2
\]
for every $\xi\in W^{1,2}$.
Since the inclusion map $W^{2,2}\to C^1$ is compact,
this is a semi-Fredholm estimate, see, e.g.,~\cite[Lemma~A.1.1]{mcduff:2004a}.
It implies finite dimension of the kernel and closedness of the image
of the operator $F_L$.
This proves Lemma~\ref{le:sF-1}.
\end{proof}

\subsubsection*{Index is conjugation invariant}

Let $\GL^+(k,\R)$ be the open subset of $\R^{k\times k}$
that consists of all real $k\times k$ matrices which are of positive
determinant, thus invertible.
Clearly $\1\in \GL^+(k,\R)$,

\begin{Lemma}\label{le:index-conj}
Consider loops of invertible matrices of positive determinant
\[
 \text{$(a)$}\
 C,G\in C^0\bigl(\SS^1,\GL^+(k,\R)\bigr)
 ,\qquad
 \text{$(b)$}\
 L,H\in W^{1,2}\bigl(\SS^1,\GL^+(k,\R)\bigr) .
\]
Then the following loops are homotopic:
\[
 C\sim GCG^{-1}
 ,\qquad
 L\sim HLH^{-1} ,
\]
and for the operators~\eqref{eq:sF-1}
the semi-Fredholm index is conjugation invariant
\[
 \INDEX F_C=\INDEX F_{GCG^{-1}}
 ,\qquad
 \INDEX F_L=\INDEX F_{HLH^{-1}} .
\]
\end{Lemma}

\begin{proof}
We prove part (a).
Part~(b) then follows by the same arguments by noting that multiplication
$W^{1,2}\times W^{1,2}\to W^{1,2}$ is continuous.

Since $C$ and $G$ take values in ${\rm GL}^+$,
we assume without loss of generality that the loops~$C$ and~$G$
are based at $\1$, in symbols $C(0)=\1=G(0)$.\footnote{If not, fix \smash{$\gamma\in C^0\bigl([0,1], GL^+(k,\R)\bigr)$}
 from $\gamma(0)=\1$ to $\gamma(1)=C(0)$.
 Define $\widebar C$ by following $\gamma$, then $C$, then
 $\gamma$ backwards (notation $\gamma^-$), all at $3$-fold speed.
 Retract $\gamma$ and $\gamma^-$ to $C_0$.}

By Lemma~\ref{le:HG-concat}, we get the based homotopies $\approx$
in Steps~1 and~3 in the calculation
\begin{equation*}
\begin{split}
 G^{-1}CG
 &\approx G^{-1}\# ( C G)
 \sim ( C G)\# G^{-1}
 \approx C GG^{-1}
 = C .
\end{split}
\end{equation*}
The (free) homotopy $\sim$ in Step 2
is the map defined for $r,t\in[0,1]$ by
\smash{$ h^r(t):= (w\# v )\bigl(\frac{r}{2}+t\bigr)$}.
It deforms $h^0=w\# v$ to $h^1=v\# w$
and moves the time $0$ point from $v(0)$ along $v$ to $v(1)$.
But the semi-Fredholm index is a homotopy invariant;
e.g.,~\cite[Section~18, Corollary~3]{Muller:2007a}. This proves
Lemma~\ref{le:index-conj}.
\end{proof}

\begin{Lemma}\label{le:HG-concat}
Let $\R^{k\times k}$ be the space of real $k\times k$ matrices.
Let $G$ and $H$ be two loops in~$\R^{k\times k}$ based at $\1=G(0)=H(0)$.
Then the pointwise matrix product is based homotopic to the
concatenation, in symbols
$ HG\approx H\# G$.
\end{Lemma}

\begin{proof}
Let $r\in[0,1]$ and $t\in[0,1]$. Defining
\[
 G^r(t)
 :=
 \begin{cases}
 G((1+r)t), &\text{$t\in\bigl[0,\frac{1}{1+r}\bigr]$,}
 \\
 \1, &\text{$t\in\bigl[\frac{1}{1+r},1\bigr]$,}
 \end{cases}
\]
then $G^1(t)$ is $G(2t)$ on \smash{$\bigl[0,\frac12\bigr]$} and $\1$ on \smash{$\bigl[\frac12,1\bigr]$}.
Defining
\[
 H^r(t)
 :=
 \begin{cases}
 \1, &\text{$t\in\bigl[0,\frac{r}{1+r}\bigr]$,}
 \\
 H((1+r)t-r), &\text{$t\in\bigl[\frac{r}{1+r},1\bigr]$,}
 \end{cases}
\]
then $H^1(t)$ is $\1$ on \smash{$\bigl[0,\frac12\bigr]$} and \smash{$H\bigl(2\bigl(t-\frac12\bigr)\bigr)$} on \smash{$\bigl[\frac12,1\bigr]$}.
Note that the homotopy defined by
$
 h^r(t)
 := H^r(t) G^r(t)
$
deforms the map
\[
 t\mapsto
 h^0(t)
 =H^0(t) G^0(t)
 =H(t) G(t)
\]
to the map
\[
 t\mapsto
 h^1(t)
 =H^1(t) G^1(t)
 = (H\# G )(t) .
\]
Moreover, the homotopy moves through based loops
\[
 h^r(0)=H^r(0)G^r(0)=\1G(0)=\1
 ,\qquad
 h^r(1)=H^r(1)G^r(1)=H(1)\1=\1 .
\]
This proves Lemma~\ref{le:HG-concat}.
\end{proof}

\subsubsection*{Fredholm index zero}

Pick a $W^{1,2}$-loop $\omega$ as in~(\ref{eq:loop-omega}).
The loop $\omega$ determines a $W^{1,2}$-loop $B=B_\omega$ which,
by~(\ref{eq:B+}), takes values in $\GL^+(2n,\R)$ and which itself, by
Lemma~\ref{le:J_B}, determines a $W^{1,2}$-loop $J_{B}$
of almost complex structures on $\R^{2n}$
compatible with $\omega$ in the sense that
$g_{J_B}:=\omega(\cdot,J_B\cdot)$
is a loop of Riemannian metrics.
Convex combination deforms
the Euclidean metric to $g_{J_B}$, namely
\begin{equation}\label{eq:hghgj7758}
 (1-r)\langle{\cdot},{\cdot}\rangle+ r g_{J_B}
 =:g^{(r)}
 =\omega\bigl(\cdot, B^{(r)}\cdot\bigr)
 ,\qquad r\in[0,1] .
\end{equation}
The identity determines a homotopy $r\mapsto B^{(r)}$
from the loop $B^{(0)}=B$ to $B^{(1)}=J_B$.
Abbreviate
\[
 J(t):=J_{B(t)}
 ,\qquad
 t\in\SS^1=\R/\Z .
\]

\begin{Lemma}\label{le:loop-exists}
There exists a loop $\Psi\in W^{1,2}\bigl(\SS^1, \GL^+(2n,\R)\bigr)$ such that
\begin{equation}\label{eq:hghgj775800}
 \Psi(t)J(0)\Psi(t)^{-1}=J(t)
 ,\qquad
 \forall t\in\SS^1 .
\end{equation}
\end{Lemma}

\begin{proof}
The proof has three steps.

\medskip\noindent
\textbf{Step~1.}
There exists a unique map
$
 \psi\colon\SS^1\to \GL(2n,\R)/\GL(n,\C)
$
such that for any $t\in\SS^1$ and any representative
$\Psi(t)\in \GL(2n,\R)$ of $\psi(t)$, we~have
\begin{equation}\label{eq:act-right}
 \Psi(t)J(0)\Psi(t)^{-1}=J(t) .
\end{equation}

\begin{proof}
Since two complex structures are conjugated,
see, e.g.,~\cite[Proposition~2.5.1]{mcduff:2004a},
for every~$t$ there exists $\Psi(t)\in \GL(2n,\R)$ such that
$\Psi(t)J(0)\Psi(t)^{-1}=J(t)$.

If \smash{$\tilde\Psi(t)$} is another element of $\GL(2n,\R)$ such that
\smash{$\tilde\Psi(t)J(0)\tilde\Psi(t)^{-1}=J(t)$}, then
\[
 \tilde\Psi(t)^{-1}\Psi(t) J(0)\bigl(\tilde\Psi(t)^{-1}\Psi(t)\bigr)^{-1}
 =J(0) .
\]
This shows that \smash{$\tilde\Psi(t)^{-1}\Psi(t)$} preserves $J(0)$ and
hence this is an element of $\GL(n,\C)$.
Therefore, the map denoted and defined by
\begin{equation}\label{eq:psi}
 \psi\colon\
 \SS^1\to\frac{\GL(2n,\R)}{\GL(n,\C)}
 ,\qquad
 t\mapsto [\Psi(t)]
\end{equation}
is well defined, namely independent of the choice of $\Psi(t)$.
In view of~(\ref{eq:act-right}),
the subgroup $\GL(n,\C)$ acts from the right on $\GL(2n,\R)$.
This proves Step~1.
\end{proof}

\noindent
\textbf{Step~2.}
The canonical map $\psi$ defined by~(\ref{eq:psi})
is of class $W^{1,2}$.

\begin{proof}
Let $t\in\SS^1$ and choose vectors $v_1,\dots,v_n\in\R^{2n}$ such that
$v_1,J(t)v_1,\dots,v_n,J(t)v_n$ is a basis of $\R^{2n}$.
Since $J$ is of class $W^{1,2}$, it is continuous.
Therefore, since linear independence is an open property, there exists
$\eps>0$ such that for any $t^\prime\in(t-\eps,t+\eps)$
the vectors $v_1,J(t^\prime)v_1,\dots,v_n,J(t^\prime)v_n$
still form a basis of $\R^{2n}$.
Choose further vectors $w_1,\dots,w_n\in\R^{2n}$
such that $w_1,J(0)w_1,\dots,w_n,J(0)w_n$ is a basis of $\R^{2n}$.
For $t^\prime\in(t-\eps,t+\eps)$, define
$\Psi(t^\prime)\in \GL(2n,\R)$ by the requirement
\[
 \Psi(t^\prime)w_j=v_j
 ,\qquad
 \Psi(t^\prime)J(0)w_j=J(t^\prime)v_j
 ,\qquad
 j=1,\dots, n .
\]
These conditions are equivalent to
\[
 \Psi(t^\prime)J(0)=J(t^\prime)\Psi(t^\prime).
\]
Moreover, since $t\mapsto J(t)$ is of class $W^{1,2}$,
it follows that the map
\[
 (t-\eps,t+\eps)\to \GL(2n,\R)
 ,\qquad
 t^\prime\mapsto\Psi(t^\prime),
\]
which locally represents $\psi$, is of class $W^{1,2}$, too.
Therefore, the map $\psi$ is locally $W^{1,2}$ and since~$\SS^1$ is
compact, it is globally $W^{1,2}$.
This proves Step~2.
\end{proof}

\noindent
\textbf{Step~3.}
The loop $\psi$ in~(\ref{eq:psi})
lifts to a $W^{1,2}$-loop $\Psi$ in
$\GL^+(2n,\R)$ based at $\1$.

\begin{proof}
As the Lie group $\GL(n,\C)$ is closed in $\GL(2n,\R)$
we have a fiber bundle
\[
 \GL(n,\C)\INTO \GL(2n,\R)\to\frac{\GL(2n,\R)}{\GL(n,\C)}
\]
see, e.g.,~\cite[Chapter~II, end of Section~13]{bredon:1993a}.
By choosing a connection on the fiber bundle,
we can lift the path $\psi$, see~(\ref{eq:psi}),
to a path $\tilde\Psi\in W^{1,2}([0,1], \GL(2n,\R))$
such that $\tilde\Psi(0)=\1$ and
\smash{$\bigl[\tilde\Psi(t)\bigr]=\psi(t)$} for every $t\in[0,1]$.
The path $\tilde\Psi$ is not necessarily a loop.
But since
$\smash{\bigl[\tilde\Psi(1)\bigr]}=\psi(1)=\psi(0)=\smash{\bigl[\tilde\Psi(0)\bigr]}$,
the initial and the end point of $\tilde\Psi$ differ by an element in $\GL(n,\C)$.
Hence there exists $\Phi_1\in \GL(n,\C)$ such that
$\tilde\Psi(1)=\tilde\Psi(0)\Phi_1=\Phi_1$.
Since $\GL(n,\C)$ is connected,
there exists a smooth path $\Phi\colon[0,1]\to \GL(n,\C)$
from $\Phi(0)=\1$ to $\Phi(1)=\Phi_1$.
The path defined for $t\in[0,1]$ by
\[
 \Psi(t):=\tilde\Psi(t)\Phi(t)^{-1}
\]
has the following properties.
Firstly, it is a loop in $\GL^+(2n,\R)$ since
\smash{$
\Psi(0):=\tilde\Psi(0)\Phi(0)^{-1}=\1
$}
and
\smash{$
\Psi(1):=\tilde\Psi(1)\Phi(1)^{-1}=\Phi(1)\Phi(1)^{-1}=\1
$}.
Secondly, it is of class $W^{1,2}$ since
$\tilde\Psi$ is of class~$W^{1,2}$ and~$\Phi$ is smooth.
Thirdly, it is a lift of $\psi$ since
\smash{$[\Psi(t)]=\bigl[\tilde\Psi(t)\Phi(t)^{-1}\bigr]=\bigl[\tilde\Psi(t)\bigr]=\psi(t)$}.
Therefore,~(\ref{eq:hghgj775800}) holds by~(\ref{eq:act-right}).
This proves Step~3.
\end{proof}

This proves Lemma~\ref{le:loop-exists}.
\end{proof}

\subsubsection*{Proof of Theorem~\ref{thm:Fredholm-F_C}}

We show Theorem~\ref{thm:Fredholm-F_C} for $k=2$, the case $k=1$
is analogous. To this end,
given a loop \smash{$\omega\in W^{1,2}\bigl(\SS^1,\Omega\bigl(\R^{2n}\bigr)\bigr)$},
we get a $W^{1,2}$-loop $B=B_\omega$ in $\GL^+(2n,\R)$.
Consider the operator
\[
 F_{B^{-1}}=B^{-1}\p_t\colon\ W^{2,2}\bigl(\SS^1,\R^{2n}\bigr)
 \to W^{1,2}\bigl(\SS^1,\R^{2n}\bigr) .
\]
By Lemma~\ref{le:sF-1}, this is a semi-Fredholm operator.
By~(\ref{eq:hghgj7758}), the loop $t\mapsto B(t)$ is homotopic in
$\GL^+(2n,\R)$ to the loop $t\mapsto J_{B(t)}$ of almost complex structures.
Therefore, the loop $B^{-1}$ is homotopic to the loop \smash{$(J_B)^{-1}=-J_B$}.

Since the semi-Fredholm index is a homotopy invariant,
see, e.g.,~\cite[Section~18, Corollary~3]{Muller:2007a},
we have
$ \INDEX F_{B^{-1}} = \INDEX F_{-J_B}$.
By Lemma~\ref{le:loop-exists}, the loop $J_B$ is conjugated in
$W^{1,2}\bigl(\SS^1,\GL^+(2n,\R)\bigr)$ to the constant loop $J(0)=J_{B(0)}$.
And therefore by Lemma~\ref{le:index-conj},
\[
 \INDEX F_{-J_B} = \INDEX F_{-J(0)}.
\]
Since $F_{-J(0)}=-J(0)\p_t$ is symmetric for the $W^{1,2}$
inner product, it holds
$
 \INDEX F_{-J(0)} =0$.
Summarizing
\[
 \INDEX F_{B^{-1}}
 =
 \INDEX F_{-J_B}
 =
 \INDEX F_{-J(0)}
 =0 .
\]
Since the index of the semi-Fredholm operator
$F_{B^{-1}}$ is zero, the operator $F_{B^{-1}}$ is actually Fredholm.
This proves Theorem~\ref{thm:Fredholm-F_C}.

\section{Heron square root iteration}

Consider the cone $\Pp$ of positive symmetric $k \times k$ matrices.
The square root map $\Pp\to\Pp$, $Q\mapsto \sqrt{Q}$,
is a diffeomorphism.
There are two ways to construct the square root map,
either using spectral calculus or the Heron iteration method.
The disadvantage of the spectral calculus approach
is that smoothness of the square root map is rather unexpected
since eigenvalues in general depend only continuously on the matrix
and the projections to the eigenspaces can even be discontinuous,
if different branches of eigenvalues meet.
Therefore, we define the square root map in this appendix with the
help of the Heron iteration.

\subsection{Real values}

\subsubsection{Heron iteration}

\begin{Lemma}[Heron method]\label{le:Heron-real}
Let $q>0$ be a positive real number.
Pick a positive real number $r_1>0$.
Define a sequence $r_n$ recursively by the requirement
\begin{equation}\label{eq:Heron-real}
 r_1>0
 ,\qquad
 r_{n+1}:=\frac12\bigl(r_n+r_n^{-1} q\bigr)
 ,\qquad
 n\in\N .
\end{equation}
Then the following is true:
\begin{itemize}\itemsep=0ex
\item[$(i)$]
 The sequence $(r_n)_{n\ge2}$ is monotone decreasing and bounded
 below by $\sqrt{q}$.
\item[$(ii)$]
 The limit $r:=\lim_{n\to\infty}r_n$ exists and $r^2=q$.
\end{itemize}
\end{Lemma}

\begin{proof}[Proof of Lemma~\ref{le:Heron-real}]
(i) We first show that
\[
 \text{(a)}\ r_n>0,
 \qquad
 \text{(b)}\ r_n\ge\sqrt{q},
\]
whenever $n\ge 2$.
While (a) is obvious from~(\ref{eq:Heron-real}), to see (b) note that
\begin{gather*}
 a:=\sqrt{\frac{r_{n-1}}{2}},\qquad
 b\stackrel{\text{(a)}}{:=}\sqrt{\frac{q}{2r_{n-1}}}
 , \quad\Rightarrow\quad
 \frac12\sqrt{q}
 =ab\le \frac{a^2+b^2}{2}
 =\frac12\biggl(\frac{r_{n-1}}{2}+\frac{q}{2r_{n-1}}\biggr)
 =\frac{r_n}{2}.
\end{gather*}
We show that the sequence is monotone decreasing.
Indeed, by (a) and (b), we get
\[
 r_n-r_{n+1}
 =r_n-\frac12\bigl(r_n+r_n^{-1} q\bigr)
 =\frac12 \frac{r_n^2-q}{r_n}\ge 0
\]
whenever $n\ge 2$.

(ii) follows from (i).
By the monotone convergence theorem,
the sequence $r_n$ has a limit $r\ge\sqrt{q}$.
On the other hand, by the recursion formula~(\ref{eq:Heron-real}) the
limit $r$ satisfies
\smash{$
 r=\frac12\bigl(r+r^{-1} q\bigr)
$},
equivalently $r^2=q$.
\end{proof}

\subsubsection{Newton--Picard iteration}

\begin{figure}[t]
 \centering
 \includegraphics
 {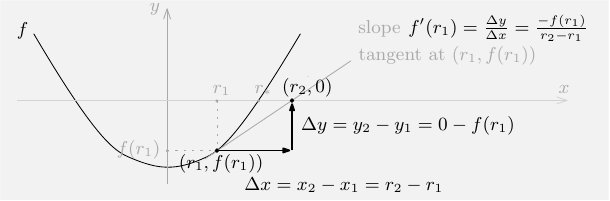}
 \caption{Newton--Picard iteration for $f(r)=r^2-q$.}
 \label{fig:fig-Newton-Picard}
\end{figure}

To determine the square root of a positive real number
$q>0$ is equivalent, to show that the function
$f(r)=r^2-q$ has a unique positive zero $r_*>0$.
For the latter, Newton--Picard iteration serves.
Choose a point $r_n>0$ and consider the tangent
line to the graph of $f$ at the point $(r_n,f(r_n))$.
If the slope $f^\prime(r_n)\not=0$ is non-zero,
the tangent line intersects the $x$-axis at a point denoted and given by
\[
 r_{n+1}=r_n-\frac{f(r_n)}{f^\prime(r_n)}
\]
as illustrated by Figure~\ref{fig:fig-Newton-Picard}.
Since $f(r_n)=r_n^2-q$ and $f^\prime(r_n)=2r_n$,
Newton--Picard iteration for the function $f(r)=r^2-q$
reproduces the Heron method, indeed
\[
 r_1>0
 ,\qquad
 r_{n+1}
 =r_n-\frac{r_n^2-q}{2r_n}
 =\frac12\biggl(r_n+\frac{q}{r_n}\biggr)
 ,\qquad n\in\N .
\]

\subsection{Matrix values}

\begin{Lemma}[Heron method]\label{le:Heron}
Let $Q\in\R^{k\times k}$ be a symmetric positive definite matrix.
Define a~sequence $R_n$ of matrices recursively by the requirement
\begin{equation}\label{eq:Heron}
 R_1:=\1
 ,\qquad
 R_{n+1}:=\frac12\bigl(R_n+R_n^{-1} Q\bigr)
 ,\qquad n\in\N .
\end{equation}
Then the following is true:
\begin{itemize} \itemsep=0ex
\item[$(i)$]
 Each matrix $R_n$ commutes with $Q$ and is symmetric positive definite.
\item[$(ii)$]
 The limit $R:=\lim_{n\to\infty}R_n$ exists and is symmetric.
\item[$(iii)$]
 $R^2=Q$ and $R$ is positive definite.
\end{itemize}
\end{Lemma}

\begin{Remark}
The symmetric positive definite matrix $R$ is called the (positive)
\textit{square root} of~$Q$ and denoted by $\sqrt{Q}:=R$.
\end{Remark}

\begin{Corollary}\label{cor:Heron}
Let $Q\in\R^{k\times k}$ be a symmetric positive definite matrix.
Then any matrix $B$ that commutes with $Q$ also commutes with $\sqrt{Q}$.
\end{Corollary}

\begin{proof}[Proof of Corollary~\ref{cor:Heron}]
Let $(R_n)_{n\in\N}$ be the sequence of positive definite symmetric
matrices defined by~(\ref{eq:Heron}).
In particular, $\lim_{n\to\infty} R_n=R=: \sqrt{Q}$.
Clearly, $B$ commutes with $R_1=\1$.
We show inductively that $B$ commutes with $R_{n+1}$ whenever $n\in\N$.
To this end, assume that $[B,R_n]=0$.
Therefore, \smash{$\bigl[B,R_n^{-1}\bigr]=0$}.
Together with the corollary hypothesis $[B,Q]=0$, we obtain that
\smash{$\bigl[B,R_n^{-1} Q\bigr]=0$}.
Hence by the recursion formula~(\ref{eq:Heron}),
we get that
\[
[B,R_{n+1}]=\frac12 \bigl([B,R_n]+\bigl[B,R_n^{-1} Q\bigr] \bigr)=0.
\]
This finishes the induction.
Hence ${[B,R]=\lim_{n\to\infty}[B,R_n]=0}$.
This proves Corollary~\ref{cor:Heron}.
\end{proof}

\begin{proof}[Proof of Lemma~\ref{le:Heron}]
(i) The proof is by induction on $n$. For $n=1$, this is true.
Suppose (i) holds for $R_n$. In particular, $R_n^{-1}$ commutes with $Q$.
Then
\begin{equation*}
\begin{split}
 R_{n+1}Q
 &=\frac12 \bigl(R_nQ+R_n^{-1} QQ \bigr)
 =\frac12 \bigl(QR_n+QR_n^{-1} Q \bigr)
 =QR_{n+1}
\end{split}
\end{equation*}
and the transpose satisfies
\[
 {R_{n+1}}^T
 \stackrel{1}{=}\frac12\bigl({R_n}^T+\bigl(QR_n^{-1}\bigr)^T\bigr)
 =\frac12\bigl(R_n+\bigl(R_n^{-1}\bigr)^TQ^T\bigr)
 \stackrel{3}{=} R_{n+1}
 .
\]
In Step 1, we used $R_n^{-1}Q=QR_n^{-1}$,
and in Step~3 that transpose of inverse is inverse of transpose.
That $R_{n+1}$ defined by~(\ref{eq:Heron})
is positive definite is true since positive definiteness
is preserved under composition of commuting symmetric matrices
and under sum of symmetric matrices.

(ii)
We follow~\cite{Higham:1986a}.
Since $Q$ is symmetric positive definite it is diagonalizable, that is,
there exists an orthogonal matrix $P$ such that
\[
 PQP^{-1}=\diag\bigl(q^{(1)},\dots,q^{(k)}\bigr)=:\Lambda,
\]
where
\smash{$q^{(1)},\dots,q^{(k)}\in\R\setminus\{0\}$} are the eigenvalues of $Q$.
Now the iterations
\[
 PR_nP^{-1}=\diag\bigl(r_n^{(1)},\dots,r_n^{(k)}\bigr)=:D_n
\]
are diagonal as well, as follows by induction on $n$.
For $n=1$, this is true since
\[
D_1:=PR_1P^{-1}=\1.
\]
Suppose it is true for $n$, then using~(\ref{eq:Heron}) in Step~2
\begin{equation*}
\begin{split}
 D_{n+1}:
 =PR_{n+1}P^{-1}
 =\frac12\bigl(PR_nP^{-1}+PR_n^{-1}P^{-1}P QP^{-1}\bigr)
 =\frac12\bigl(D_n+D_n^{-1} \Lambda\bigr)
\end{split}
\end{equation*}
is indeed a diagonal matrix.
But now each diagonal position corresponds to a real-valued Heron
iteration~(\ref{eq:Heron-real}). This proves that the limit $R$ exists.

Symmetry: Given $\xi,\eta\in\R^k$, then
$\INNER{R_n\xi}{\eta}=\INNER{\xi}{R_n\eta}$ by symmetry of $R_n$.
In the limit, as~${n\to\infty}$, we get
$\INNER{R\xi}{\eta}=\INNER{\xi}{R\eta}$.

(iii)~By~(ii), the limit $R$ exists and by the recursion
formula~(\ref{eq:Heron}) it satisfies
$R=\frac12\bigl(R+R^{-1} Q\bigr)$,
equivalently $R^2=Q$.

Positive definite:
Since the $R_n$ are symmetric positive definite,
their eigenvalues are real and~${\ge 0}$.
Since $R=\lim_{n\to\infty} R_n$ is symmetric,
its eigenvalues are real. Since eigenvalues depend continuously
on the matrix the eigenvalues of $R$ are $\ge 0$.
Let $\rho$ be the smallest eigenvalue of~$R$.
Then $\rho\ge 0$ and $\rho^2$ is eigenvalue of $R^2=Q$.
Since $Q$ is positive definite, all eigenvalues are
strictly positive, in particular~$\rho^2>0$, hence $\rho>0$.
Thus~$R$ is positive definite.

This concludes the proof of Lemma~\ref{le:Heron}.
\end{proof}

If in Corollary~\ref{cor:Heron}
one assumes in addition symmetric and positive definite for the matrix
$B$, then the proof reduces essentially to the existence of a basis
composed of common eigenvectors of $B$ and $Q$.

\begin{Lemma}\label{le:commute-square}
Let $Q\in\R^{k\times k}$ be a symmetric positive definite matrix.
Then a symmetric positive definite matrix $B$ that commutes with~$Q$
commutes with $\sqrt{Q}$.
\end{Lemma}

\begin{proof}
Linear algebra tells the following:
Two symmetric positive definite $k\times k$ matrices, in the case at
hand $Q$ and $B$, commute if and only if there exists an orthonormal basis
$\Xx=\{\xi_1,\dots,\xi_k\}$ whose elements are eigenvectors of both
matrices, say $Q\xi_i=\rho_i\xi_i$ and $B\xi_i=\lambda_i\xi_i$ for
some $\rho_i,\lambda_i>0$ and $i=1,\dots,k$.
The positive square root of $Q$, notation $\sqrt{Q}=:R$, is defined
by $R\xi_i:=\sqrt{\rho_i}\xi_i$ for $i=1,\dots,k$.
It is an exercise to check that $R$ is symmetric positive definite;
here pairwise orthogonality of the $\xi_i$ enters.
But the ON basis $\Xx$ is composed of eigenvectors of both~$R$
and~$B$, hence both matrices commute by the linear algebra
cited in the beginning of this proof.

This proves Lemma~\ref{le:commute-square}.
\end{proof}

\subsection*{Acknowledgements}
UF~acknowledges support by DFG grant FR~2637/5-1.

\pdfbookmark[1]{References}{ref}
\LastPageEnding

\end{document}